\documentstyle[amssymb,amsfonts]{amsart}

\newenvironment{proof}{\noindent {\bf Proof} }{\endprf\par}
\def \endprf{\hfill  {\vrule height6pt width6pt depth0pt}\medskip}
\def\emph#1{{\it #1}}
\def\textbf#1{{\bf #1}}

\newcommand{\bea}{\begin{eqnarray}}
\newcommand{\eea}{\end{eqnarray}}
\def\beaa{\begin{eqnarray*}}
\def\eeaa{\end{eqnarray*}}

\def\ba{\begin{array}}
\def\ea{\end{array}}
\def\be#1{\begin{equation} \label{#1}}

\newcommand{\nn}{\nonumber}
\def\medn{\medskip\noindent}
\def\rrrr{{\Bbb R}}
\def\rr{{\bf R}}
\def\rb{{\bf w}}

\def\nn{\nonumber}
\def\stu{{S_{t,u}}}
\def\sttu{{S_{t',u}}}
\def\stau{{S_{\tau,u}}}

\def\logp{\log^+}

\def\nabc{\overset{\circ}{\dd}}
\def\Lie{{\cal L}}
\def\tr{\mbox{tr}}

\def \Gmu{G^{(\mu)}}
\def\MH{{\cal M}(\pr H)}
\def\gg{{\bf g}}
\def\mm{{\bf m}}

\newcommand\und{\underline}
\newcommand{\half}{\frac 12}

\def\f12{\frac 12}
\def\half{\frac 12}

\def\bk{\bar{k}}

\newcommand{\divv}{\mbox{div}\mkern-19mu /\,\,\,\,}
\newcommand{\curll}{\mbox{curl}\mkern-19mu /\,\,\,\,}
\newcommand{\muv}{\mu\mkern-9mu /}

\def\a{\alpha}
\def\b{\beta}
\def\ga{\gamma}
\def\Ga{\Gamma}
\def\de{\delta}

\def\ep{\epsilon}
\def\eps{\epsilon}
\def\la{\lambda}
\def\La{\Lambda}
\def\si{\sigma}
\def\Si{\Sigma}

\def\Om{\Omega}

\def\th{\theta}

\def\nab{\nabla}

\newcommand{\trchb}{\tr \chib}
\newcommand{\chih}{\hat{\chi}}
\newcommand{\chib}{\underline{\chi}}
\newcommand{\xib}{\underline{\xi}}
\newcommand{\etab}{\underline{\eta}}
\newcommand{\chibh}{\underline{\hat{\chi}}}
\def\f14{\frac{1}{4}}
\def\dd{{\bf D}}

\newcommand{\les}{\lesssim}

\newcommand{\piT}{\,^{(T)}\pi}

\def\il{\int}

\def\Lb{\underline{L}}

\def\ric{\mbox{ \bf Ric}}

\def\pr{\partial}
\def\f12{\frac 1 2}
\renewcommand{\chib}{\underline{\chi}}
\def\chih{\hat{\chi}}
\def\trch{\mbox{tr}\chi}
\newcommand{\ddd}{{\bf \mathcal D} \mkern-13mu /\,}
\newcommand{\nabb}{\mbox{$\nabla \mkern-13mu /$\,}}
\def\ah{\hat{\alpha}}

\begin{document}
\theoremstyle{plain}
  \newtheorem{theorem}[subsection]{Theorem}
  \newtheorem{conjecture}[subsection]{Conjecture}
  \newtheorem{proposition}[subsection]{Proposition}
  \newtheorem{lemma}[subsection]{Lemma}
  \newtheorem{corollary}[subsection]{Corollary}

\theoremstyle{remark}
  \newtheorem{remark}[subsection]{Remark}
  \newtheorem{remarks}[subsection]{Remarks}

\theoremstyle{definition}
  \newtheorem{definition}[subsection]{Definition}

\include{psfig}

\title[rough Einstein metrics ]{The causal structure of microlocalized rough
Einstein metrics}
\author{Sergiu Klainerman}
\address{Department of Mathematics, Princeton University,
 Princeton NJ 08544}
\email{ seri@@math.princeton.edu}

\author{Igor Rodnianski}
\address{Department of Mathematics, Princeton University, 
Princeton NJ 08544}
\email{ irod@@math.princeton.edu}
\subjclass{35J10}
\vspace{-0.3in}
\begin{abstract}
This is the second   in a series of three  papers in which 
we initiate the study of very rough  solutions  to the initial value problem
for the Einstein vacuum equations expressed relative
 to  wave coordinates. By very rough we mean solutions
which cannot be constructed by the classical techniques of energy estimates 
and Sobolev inequalities. In this paper we develop
the  geometric analysis of  the Eikonal equation for   microlocalized
 rough Einstein metrics. This is a  crucial step in the derivation
of the decay estimates needed in the first paper.

\end{abstract}
\maketitle
\section{Introduction} 
This is the second in a series of three papers in which we 
initiate the study
of \textit{ very rough} solutions of the Einstein vacuum equations.
 By very  rough we mean solutions which can not be dealt with
by the classical techniques of energy estimates and Sobolev inequalities.
 In fact in this work we  develop and  take   advantage of 
Strichartz type estimates. The result, stated in our first paper \cite{Einst1},
is in fact optimal with respect to the full  potential of such 
estimates\footnote{To go beyond our result will require the
 development of bilinear techniques for the Einstein equations, 
see the discussion in the introduction to \cite{Einst1}.}.
We recall below our main 
result:

\begin{theorem}[Main Theorem]  Let $\gg$ be  a classical solution\footnote{We denote 
by $R_{\a\b} $ the Ricci curvature of $\gg$.} 
 of the  Einstein  equations 
\be{I1}
\rr_{\a\b}(\gg)=0
\end{equation}
expressed\footnote{In wave coordinates
the Einstein equations take the reduced  form 
$$
\gg^{\a\b}\pr_\a\pr_\b  \gg_{\mu\nu}=N_{\mu\nu}(\gg,\pr \gg)
$$
with $N$ quadratic in the first derivatives
 $\pr\gg$ of the metric.} relative to  wave coordinates $x^\a$,
\be{I2}
\square_{\gg} x^\a =\frac{1}{|\gg|}\pr_\mu(\gg^{\mu\nu}|\gg|\pr_\nu)x^\a= 0.
\end{equation}
We assume that on  the initial spacelike  hyperplane  $\Sigma$ 
 given by $t=x^0=0$,
$$
\nabla \gg_{\a\b}(0)\in  H^{s-1}(\Sigma)\,,\quad  \pr_t \gg_{\a\b}(0)\in 
H^{s-1}(\Sigma)
$$
with $\nabla$ denoting the gradient with respect to the
 space coordinates $x^i$, $i=1,2,3$ and  $H^s$ the standard 
Sobolev spaces. We also assume that $\gg_{\a\b}(0)$ is
a continuous Lorentz metric and 
$
\sup_{|x|=r}|\gg_{\a\b}(0)-\mm_{\a\b}|\longrightarrow 0\quad  \mbox{as}  \quad
r\longrightarrow \infty,
$
where $|x|=(\sum_{i=1}^3 |x^i|^2)^{\frac{1}{2}}$ and $\mm_{\a\b}$
the Minkowski metric.

 We show\footnote{We assume however
 that $T$ stays sufficiently small, e.g. $T\le 1$.  This a purely technical
 assumption which one should be able to remove.} that  the
time
$T$ of existence depends in fact only on the size of the norm $\|\pr
\gg_{\mu\nu}(0)\|_{H^{s-1}}$,
 for any fixed   $s>2$. 
\end{theorem}

In \cite{Einst1} we have given a detailed proof of the Theorem
by relying heavily on a result, we have called the Asymptotic Theorem,
concerning the geometric properties of the causal structure of appropriately
microlocalized rough Einstein metrics. This result, which is  the focus
of  this paper, is of independent interest as it requires the
development of   new geometric and analytic methods to deal with 
characteristic surfaces of the Einstein metrics.

More precisely we study the solutions, called optical functions,  of the Eikonal equation
\be{eik}
H_{(\la)}^{\a\b}\pr_\a u\pr_\b u=0,
\end{equation}
associated to the family of  regularized Lorentz metrics $H_{(\la)}$, $\la\in 2^{\Bbb N}$,
defined, starting with an $H^{2+\eps}$ Einstein metric $\gg$, by the formula
\be{1.100}
H_{(\la)}= P_{<\la} \gg(\la^{-1} t ,\la^{-1} x).
\end{equation}
where\footnote{More precisely, for a given function of the spatial variables $x=x^1,x^2,x^3$,
 the Littlewood Paley projection $P_{<\la}f=
\sum_{\mu< \f12 \la} P_\mu f$, \,\, $P_\mu f={\cal
F}^{-1}\big(\chi(\mu^{-1}\xi)\hat{f}(\xi)\big)$ with $\chi$ supported in the unit  dyadic
region $\half \le |\xi|\le 2$.} $P_{<\la}$ is an operator which cuts 
off all the frequencies above\footnote{The definition of the projector $P_{<\la}$ in 
\cite{Einst1} was slightly different from the one we are using in this paper.
 There $P_{<\la}$ removed all the frequencies above $2^{-M_0}\la$ for some sufficiently
large constant $M_0$. It is clear that a simple rescaling can remedy this
discrepancy.} 
 $\la$.
The importance of  the eikonal equation \eqref{eik} in the study of solutions
 to wave equations on a background Lorentz metric $H$ is well known.
It is mainly used,  in the geometric optics approximation, to construct
parametrices associated to the corresponding linear operator $\square_H$.
In particular it has played a fundamental  role in the recent works 
of  Smith\cite{Sm},  Bahouri-Chemin \cite{Ba-Ch1}, \cite{Ba-Ch2} 
 and Tataru \cite{Ta1},\cite{Ta2} concerning rough solutions to linear and 
nonlinear wave equations. Their work relies indeed on parametrices defined
with the help of specific families of optical functions corresponding 
to null hyperplanes. In \cite{Kl}, \cite{Kl-Ro}, and also \cite{Einst1}  which do not rely
on specific parametrices, a special optical function, corresponding
to null cones with vertices on a timelike geodesic,  was used to construct
an almost conformal Killing vectorfield.  

 The main message of our paper is that optical functions 
 associated to Einstein metrics, or microlocalized  versions
of them, have better  properties. This fact was already
recognized in \cite{Ch-Kl} where the construction of an optical function
normalized at infinity played a crucial role in the proof of the
 global nonlinear stability  of the Minkowski space. A similar construction,
based on two optical functions, can be found in \cite{Kl-Ni}. Here,
 we take the use of the  special structure of the Einstein equations
one step further by deriving unexpected regularity properties
of optical functions which are essential in the proof
of the Main Theorem.   It was well known (see \cite{Ch-Kl}, \cite{Kl},
\cite{Kl-Ro}) that the use of Codazzi equations combined with the Raychaudhuri 
equation for the $\trch$, the trace of null second fundamental form $\chi$,
leads to the improved estimate
for the first angular derivatives of the traceless part of $\chi$. 
A similar observation holds for another null component of the Hessian of
the optical function, $\eta$. The role of the Raychaudhuri equation is taken
by the transport equation  the 
``mass aspect function'' $\mu$. 

\noindent In this paper we show, using the structure of the
curvature terms in the main equations, how to derive improved regularity estimates
for the undifferentiated quantities $\chih$ and $\eta$.  
In particular, in the case of the estimates for $\eta$ we are lead to 
introduce a new non local
quantity $\muv$ tied to $\mu$ via a Hodge system.

The properties of the optical function are given in details 
in the statement of the Asymptotics Theorem. We shall give
a precise statement of it in section 2 after we introduce a few 
essential definitions. 

The paper is organized as follows:
\begin{itemize}
\item In section \ref{se:geompr} we construct an optical function $u$,
constant on null cones with vertices on a fixed timelike geodesic, and
 describe our  basic geometric entities
associated  to it.  We define the surfaces $S_{t,u}$, the canonical  
null pair  $L, \Lb$ and the associated  Ricci coefficients. This allows
us   to give a  precise statement of our main
result, the Asymptotic Theorem \ref{ITrih}. 
\item In section \ref{se:nullstr} we derive the structure equations for
the Ricci coefficients. These equations are  a coupled system of the transport
and Codazzi equations and are fundamental for the proof of
theorem\ref{ITrih}.
\item In section \ref{se:specialstr} we obtain some crucial properties
of the components of the Riemann curvature tensor $\rr_{\a\b\ga\de}$.
 
\item The remaining sections are occupied with the proof of the Asymptotic 
Theorem. We give a  detailed description of their content 
and strategy of the proof  in
section \ref{se:strategy}.
\end{itemize}
The paper is essentially self-contained. From the first paper in this series
\cite{Einst1} we only need the result of proposition 2.4 (Background Estimates)
which in any case can be easily derived from the the metric hypothesis 
\eqref{bootstrap}, the Ricci condition \eqref{I1}, and the definition \eqref{1.100}.
We do however rely on the following results which will be proved in a forthcoming 
paper \cite{Einst3}:
\begin{itemize}
\item Isoperimetric and trace inequalities, see proposition \ref{Triso}
\item  Calderon-Zygmund type estimates, see proposition \ref{CZestimates}
\item Theorem \ref{Ricci4}
\end{itemize}

We recall our metric hypothesis( referred in  \cite{Einst1}, section 2
 as the  bootstrap  hypothesis)  on the components
of
$\gg$
 relative to our wave  coordinates $x^\a$,

\medn
{\bf Metric Hypothesis:}
\be{bootstrap} 
\|\pr\gg\|_{L^\infty_{[0,T]} H^{1+\ga}} +\|\pr \gg\|_{L^2_{[0,T]}
L^\infty_x}\le B_0,
\end{equation}
for some fixed $\ga>0$.
\section{ Geometric preliminaries  }
\label{se:geompr}
We start by recalling the basic geometric constructions
associated with a Lorentz metric $H=H_{(\la)}$.

 Recall, see \cite{Einst1} section 2,  that  the parameters of the $\Si_t$
foliation are given by $n, v$, the induced  metric $h$ and the second
fundamental form $k_{ij}$, according to the decomposition,
\be{2.3}
H=-n^2 dt^2+h_{ij}(dx^i+v^i dt)\otimes(dx^j+v^j dt),
 \end{equation} 
with $h_{ij}$ the induced Riemannian  metric on $\Si_t$,
$n$ the lapse and $v=v^i\pr_i$ the shift of $H$. Denoting by
$T$ the unit, future oriented, normal to $\Si_t$ and $k$
the second fundamental form $k_{ij}=-<\dd_iT, \pr_j>$ we 
find,
\begin{align}
 &\pr_t=n T+v,\qquad <\pr_t, v>=0\nn \\
 &k_{ij}= -\frac{1}{2}\Lie_T H_{\,ij}=-12 n^{-1}(\pr_t
 h_{ij}-\Lie_v h_{\,ij})\label{lit}
\end{align}
with $\Lie_X$ denoting the Lie derivative with respect to the 
vectorfield $X$.
 We also have 
the following, see \cite{Einst1}  sections 2, 8:
\be{4.1}
c|\xi|^2\le h_{ij} \xi^i\xi^j\le c^{-1}|\xi|^2, \qquad c \le n^2-|v|_h^2
\end{equation}
for some $c>0$.
Also
$
n, |v|\les 1.
$

The time axis is
 defined as the integral curve of the forward  unit normal 
$T$ to the hypersurfaces $\Si_t$.  The point $\Ga_t$ is
the intersection between $\Ga$ and $\Si_t$.
\begin{definition} The  optical function $u$ is 
an outgoing  solution of the Eikonal equation
\be{it1}
H^{\a\b}\pr_\a u \pr_\b u=0
\end{equation}
 with initial conditions
$u(\Ga_t)=t$ on the time axis. 
\end{definition}
  The level
surfaces of u, denoted $C_u$ are outgoing  null
 cones with vertices
on the time axis. Clearly,
\be{it1'}
T(u)=|\nab u|_h
\end{equation}
where $h$ is the induced metric on $\Si_t$, 
$|\nab u|_h^2=\sum_{i=1}^3|e_i(u)|^2$ relative to an orthonormal
frame $e_i$ on $\Si_t$. 

We denote by $S_{t,u}$  the surfaces of intersection
between $\Si_t$ and  $C_u$. They play a fundamental role
in our discussion.
\begin{definition}[\textit{  Canonical null pair}]
 \be{it2}
L=bL'=T+N, \qquad \Lb=2T-L=T-N
\end{equation}
Here $L'=-H^{\a\b}\pr_\b u \pr_\a$ is the geodesic null generator 
 of $C_u$,  $b$  is  the  \textit{lapse of the null foliation}(or shortly null lapse)
\be{it3}
b^{-1}=-<L', T>=T(u),
\end{equation} 
and $N$ the  exterior unit normal, along $\Si_t$, to the surfaces
 $S_{t,u}$. 
\end{definition}
\begin{definition} Null frame , $e_1,e_2,e_3,e_4$
\end{definition}
\begin{definition}[\textit{Ricci coefficients}]

 Let   $e_3=\Lb$, $e_4=L$ be our canonical null
pair  and 
$(e_A)_{A=1,2}$ an arbitrary  orthonormal frame\footnote{
$e_1,e_2,e_3,e_4$ forms a  null
 frame. This can always be defined locally, in a neighborhood
of a point. } on
$S_{t,u}$.  The following tensors on  $S_{t,u}$ 
\begin{alignat}{2}
&\chi_{AB}=<\dd_A e_4,e_B>, &\quad 
&\chib_{AB}=<\dd_A e_3,e_B>,\label{chi}\nn\\
&\eta_{A}=\half <\dd_{3} e_4,e_A>,&\quad
&\etab_{A}=\half <\dd_{4} e_3,e_A>,\label{eta}\,\label{xi}\\
&\xib_{A}=\half <\dd_{3} e_3,e_A>.\nn
\end{alignat}
are called the Ricci coefficients associated to our canonical null pair.

We decompose $\chi$ and $\chib$ into
their  trace and traceless components.
\begin{alignat}{2}
&\trch = H^{AB}\chi_{AB},&\quad &\trchb = H^{AB}\chib_{AB},
\label{trchi}\\
&\chih_{AB}=\chi_{AB}-\half \trch H_{AB},&\quad 
&\chibh_{AB}=\chib_{AB}-\half \trchb H_{AB},
\label{chih} 
\end{alignat}
\end{definition}

We define $s$ to be the affine parameter
of $L$, i.e. $L(s)=1$ and $s=0$ on the time axis $\Ga_t$. In \cite{Kl-Ro},
where $n=1$ we had $s=t-u$. Such a simple relation does not hold in our case,
we have instead, along any fixed $C_u$,
\be{sandt}
\frac{dt}{ds}=n^{-1}
\end{equation}
We shall also introduce the area $A(t,u)$ of the 2-surface $S(t,u)$
and the radius $r(t,u)$ defined by
\be{defr}
A=4\pi r^2 
\end{equation}
Along a given $C_u$ we have\footnote{This follows by writing
the metric on $S_{t,u}$ in the form $\ga_{AB}(s(t,\th), \th) d\th^a d\th^B$,
 relative to angular coordinates $\th^1, \th^2$, and its area $A(t,u)=\int \sqrt{\ga}
d\th^1\wedge d\th^2$. Thus, that $\frac{d}{dt}A=\int
\f12\ga^{AB}\frac{d}{dt}\ga_{AB}\sqrt{\ga} d\th^1\wedge d\th^2$.
  On the other hand $\frac{d}{ds} \ga_{AB}=2\chi_{AB}$  and $\frac{ds}{dt}=n$.}
$$\frac{\pr A}{\pr t}=\int_S n \trch.$$
Therefore, along $C_u$,
\be{randt}
\frac{dr}{dt}=\frac{r}{2}\overline{n\trch}
\end{equation}
where, given a function $f$ we denote by $\bar{f}(t,u)$ its average
on $S_{t,u}$. Thus  $$\bar{f}(t,u)=
\frac{1}{4\pi r^2}\int_{S_{t,u}} f.
$$
The following  \textit{Ricci equations} can also
be easily derived see \cite{Kl-Ro}. They express  the covariant 
derivatives $\dd$ of the null frame $(e_A)_{A=1,2}, e_3, e_4$ relative to itself.
\begin{alignat}{2}
&\dd_A e_4=\chi_{AB} e_B - k_{AN} e_4, &\quad 
&\dd_A e_3=\chib_{AB} e_B + k_{AN} e_3,\nn\\
&\dd_{4} e_4 = -\bk_{NN} e_4, &\quad 
&\dd_{4} e_3= 2\etab_{A} e_A + \bk_{NN} e_3, 
  \label{ricciform} \\
 &\dd_{3} e_4 = 2\eta_{A}e_A +
\bk_{NN} e_4,&\quad &\dd_{3} e_3 = 2\xib_{A}e_A -
\bk_{NN} e_3,\nn\\ &\dd_{4} e_A = \ddd_{4} e_A +
\etab_{A} e_4,&\quad &\dd_{3} e_A = \ddd_{3} e_A+\eta_A e_3
+ \xib_A e_4,
\nn\\
&\dd_{B}e_A = \nabb_B e_A +\half \chi_{AB}\, e_3 +\half
 \chib_{AB}\, e_4\nn
\end{alignat}
where, $\ddd_{ 3}$, $\ddd_{ 4}$ denote the 
projection on $S_{t,u}$ of $\dd_3$ and $\dd_4$, $\nabb$
denotes the induced covariant derivative on $S_{t,u}$
and, for every vector $X$ tangent to $\Si_t$, 
\be{newk}
\bk_{NX}=k_{NX}-n^{-1} \nab_X n
\end{equation}
Thus $\bk_{NN}=k_{NN}-n^{-1}N(n)$ and  $\bk_{AN}=k_{AN}-n^{_1}\nab_A n$.
Also,
\begin{align}
&\chib_{AB}=-\chi_{AB}-2k_{AB},\nn\\
&\etab_A = -\bk_{AN},\label{etab}\\
&\xib_A = k_{AN}+n^{-1} \nab_A n-\eta_A.\nn
\end{align}
and,
\begin{equation}
\eta_A = b^{-1}\nabb_A b + k_{AN}.
\label{etaa}
\end{equation}
The formulas \eqref{ricciform},
\eqref{etab} and
\eqref{etaa} can be checked
in precisely the same manner as 
(2.45--2.53) in \cite{Kl-Ro}. The only
difference occur because $\dd_TT$ does not
longer vanishes. We have in fact, relative
to any orthonormal frame $e_i$ on $\Si_t$,
\be{3.2}
\dd_TT=n^{-1}e_i(n) e_i
\end{equation}
To check \eqref{3.2} observe that we can introduce 
new local coordinates $\bar{x}^i=\bar{x}^i(t,x)$ on $\Si_t$
  which preserve  the lapse $n$  while  making  the shift $V$
 to vanish identically. Thus $\pr_t=nT$ and therefore,
for an arbitrary  vectorfield $X$ tangent to  $\Si_t$, we easily
calculate,
$<\dd_T T, X>=n^{-2} X^i<\dd_{\pr_t} \pr_t, \pr_i>=
-n^{-2}X^i< \pr_t, \dd_{\pr_t}\pr_i>
=-n^{-2}X^i<\pr_t,\dd_{\pr_i}\pr_t> 
=-n^{-2} X^i\f12 \pr_i<\pr_t,\pr_t>= n^{-2} X^i\f12
\pr_i(n^2)=n^{-1}X(n)$.

Equations \eqref{etab} indicate that the only independent geometric 
quantities, besides $n$, $v$ and $k$ are
 $\trch, \chih, \eta$.  
We now state the main result of our paper giving the precise description
of the Ricci coefficients. 
\begin{theorem}
Let $\gg$ be an Einstein metric obeying the Metric  Hypothesis \eqref{bootstrap}
and $H=H_{(\la)}$ be the family of the regularized Lorentz metrics defined
according to \eqref{1.100}. Fix a sufficiently large value of the dyadic parameter 
$\la$ and consider, corresponding to $H=H_{(\la)}$, the  
optical function $u$ defined above. Let ${\cal I}^+_0$ be the future domain
of the origin on $\Si_0$.  

Then for any $\eps_0>0$, such that $5\eps_0<\ga$ with $\ga$ from \eqref{bootstrap},
we can extend  the optical function $u$ throughout the region $ {\cal I}_0^+\cap( [0, 
\la^{1-8\eps_0}]\times
\rrrr^3) $ and show that 
in that region  the Ricci coefficients
$\trch $, $\chih$,
and $\eta$
satisfy the following estimates:
\begin{equation}
\begin{align}
&\|\trch -\frac 2r\|_{L^2_t L^\infty_x}+ \|\chih\|_{L^2_t L^\infty_x} +
 \|\eta\|_{L^2_t L^\infty_x}\les \la^{-\frac 12-3\eps_0},\label{itrih1}\\
&\|\trch -\frac 2r\|_{L^q(\stu)} + \|\chih\|_{L^q(\stu)} +
\|\eta \|_{L^q(\stu)}\les \la^{-3\eps_0}.
\label{trih222}
\end{align}
\end{equation}
In the estimate \eqref{trih1} the  function $\frac 2r$ can be replaced 
with $\frac 2{n(t-u)}$. 
In addition, in the exterior region $r\ge t/2$, 
\begin{equation}
\begin{split}
&\|\trch -\frac 2s\|_{L^\infty(\stu)}\les t^{-1} \la^{-4\eps_0},\qquad
\|\chih\|_{L^\infty(\stu)}\les t^{-1} \la^{-\eps_0} + \|\pr H(t)\|_{L^\infty_x},\\
&\|\eta\|_{L^\infty(\stu)}\les \la^{-1} + \la^{-\eps_0} t^{-1} 
+  \la^{\eps}\|\pr H(t)\|_{L^\infty_x}.
\end{split}
\label{itrih3}
\end{equation}
where the last estimate holds for an arbitrary positive $\eps$, $\eps<\eps_0$. 
We also have the following estimates for the derivatives of $\trch$:
\begin{align}
& \|\sup_{r\ge \frac t2}
\|\Lb(\trch -\frac 2{r})\|_{L^2(\stu)}\|_{L^1_t} + \|\sup_{r\ge \frac t2}
\|\Lb(\trch - \frac 2{n(t-u)})\|_{L^2(\stu)}\|_{L^1_t}\le \la^{-3\eps_0 },
\label{i2trih6}\\
&\|\sup_{r\ge \frac t2}\|\nabb \trch \|_{L^2(\stu)}\|_{L^1_t} +
\|\sup_{r\ge \frac t2}\|\nabb \big (\trch -\frac 2{n(t-u)}\big)
\|_{L^2(\stu)}\|_{L^1_t}\le \la^{-3\eps_0 }
\label{i2trih8}
\end{align}
In addition we also have  weak  estimates of the form,
\be{lastz}
\sup_{u\le \frac t2}\|(\nabb, \Lb)\big( \trch -\frac
2{n(t-u)}\big)\|_{L^\infty(\stu)}\les \la^C
\end{equation}
for some large value of $C$.
The inequalities $\les$ indicate that the bounds hold with some universal 
constants including the constant $B_0$ from \eqref{bootstrap}.
\label{ITrih}
\end{theorem}

\section{Null structure equations}
\label{se:nullstr}
In the proof of theorem\ref{ITrih} we rely on the system of equations
satisfied by the by the Ricci coefficients $\chi$, $\eta$.
Below we  write down our main structure equations. Their
derivation proceeds in exactly the same way as in 
\cite{Kl-Ro}( see propositions 2.2 and 2.3)
from  the  formulas \eqref{ricciform} above. 
\begin{proposition}
The components   $\tr\chi,\, \chih,\,\eta$ and the lapse
$b$ verify the following equations\footnote{which can be interpreted as
transport equations along the  null geodesics 
generated by $L$.}:
\begin{align}
&L(b) = - b\, \bk_{NN}, \label{D4a}\\
&L(\trch) + \half (\trch)^2 = - |\chih|^2 -\bk_{NN} \trch - 
\rr_{44},
\label{D4trchi}\\
&\ddd_{4} \chih_{AB} + \half \trch \chih_{AB} = 
-\bk_{NN} \chih_{AB} - \ah_{AB},
\label{D4chih}\\
&\ddd_{4}\eta_A + \half (\trch)\eta_A = -(k_{BN} +\eta_B)\chih_{AB} -
\half \trch k_{AN} - \half \beta_A,\label{D4eta}.
\end{align}
\label{D4}
Here $\ah_{AB}= \rr_{4A4B} - \f12 \rr_{44} \de_{AB}$ and $\b_A=\rr_{4A34}$. 
Also, setting, 
\be{eqmu}
\mu=\Lb(\trch)-\f12 (\trch)^2 - 
\big (k_{NN}+n^{-1}\nab_N n\big )\trch
\end{equation}
we find \begin{equation}
\begin{split}
L(\mu) + \trch \mu &=2(\etab_A-\eta_A)\nabb_A(\trch)
-2\chih_{AB}\Bigl (2\nabb_A \eta_B +  2\eta_A\eta_B   \\&+ \bk_{NN}\chih_{AB} 
+\trch\chih_{AB} +\chih_{AC}\chih_{CB} +2k_{AC}\chi_{CB}+ \rr_{B{43}A}\Bigr )
\\&
-\Lb(\rr_{44})+  (2k_{NN} - 4 n^{-1}\nab_N n)
)\big (\half (\trch)^2 - |\chih|^2 -
\bk_{NN} \trch -\rr_{44}\big ) \\ & +4\bk_{NN}^2\trch 
+(\trch+4\bk_{NN})(|\chih|^2+\rr_{44}) \\ &- \trch 
\bigg (2 (k_{AN} - \eta_A) n^{-1} \nab_A n - 2 |n^{-1} N(n)|^2 +
\rr_{4343} + 2 k_{Nm} k^{m}_N\bigg)
\end{split}
\label{D4tmu}
\end{equation}
\label{proptransp}
\end{proposition}
\begin{remark}
Equation \eqref{D4trchi} is known as the Raychaudhuri equation in the relativity
literature, see e.g. \cite{Ha-El}.
\end{remark}
\begin{remark} Observe that our definition of $\mu$ differs from that
in \cite{Kl-Ro}. Indeed there we had, instead of $\mu$,
$$\tilde\mu=\Lb(\trch)-\f12 (\trch)^2-3\bk_{NN}\trch$$
\end{remark}
and the corresponding transport equation:
\begin{equation}
\begin{split}
L( \tilde\mu) + \trch \tilde\mu &=2(\etab_A-\eta_A)\nabb_A(\trch)
-2\chih_{AB}\Bigl (2\nabb_A \eta_B +  2\eta_A\eta_B   \\&+ \bk_{NN}\chih_{AB} 
+\trch\chih_{AB} +\chih_{AC}\chih_{CB} +2k_{AC}\chi_{CB}+ \rr_{B{43}A}\Bigr )
\\&
-\Lb(\rr_{44})-\Lb(\bk_{NN})\trch -3L(\bk_{NN})\trch +4\bk_{NN}^2\trch 
\\&+(\trch+4\bk_{NN})(|\chih|^2+\rr_{44})
\end{split}
\label{D4mu}
\end{equation}
We obtain \eqref{D4tmu} from \eqref{D4mu} as follows:
The second fundamental form $k$ verifies the equation(
see formula (1.0.3a) in \cite{Ch-Kl}),
$$
{\mathcal L}_{nT} k_{ij} = -\nab_i \nab_j n + n(\rr_{iTjT}-k_{im} k^{m}_j).
$$
In particular,
$$
{\mathcal L}_{nT} k_{NN} = -\nab^2_N n + n(\rr_{NTNT}-k_{Nm} k^{m}_N).
$$ 
Exploiting the definition of the Lie derivative ${\mathcal L}_{nT}$, we obtain
$$
T(k_{NN}) + 2  k(\nab_N T, N) = -n^{-1} \nab^2_N n + (\rr_{NTNT}-k_{Nm} k^{m}_N).
$$
It then follows that
$$
\half \Lb (k_{NN}) + \half L(k_{NN}) - 2 (k_{NN})^2 - 2(k_{AN})^2 =  -n^{-1} \nab^2_N n + 
(\rr_{NTNT}-k_{Nm} k^{m}_N) 
$$
Therefore,
$$
\aligned
\half \Lb (k_{NN}) - \half \Lb \big (n^{-1} N (n)\big ) & =
- \half L(k_{NN}) - \half L \big (n^{-1} N (n)\big ) + 
(\rr_{NTNT}+ k_{Nm} k^{m}_N)\\ &+ n^{-1}(\nab_N N) n  - n^{-2} |N(n)|^2
\endaligned
$$
Recall that $\bk_{NN}= k_{NN} -n^{-1} N (n)$ and 
$<\nab_N N, e_A> = k_{AN} -\eta_A$. Thus 
\beaa
\Lb (\bk_{NN})& =& - L \big( k_{NN} + n^{-1} N (n)\big ) + 
2 (k_{AN} - \eta_A) n^{-1} \nab_A n\\
& -& 2 |n^{-1} N(n)|^2 +
\rr_{4343} + 2 k_{Nm} k^{m}_N.
\eeaa 
Therefore taking  $ \mu =\Lb(\trch)-\f12 (\trch)^2- (k_{NN} + n^{-1} N (n))\trch$
we derive the desired transport equation \eqref{D4tmu}.

\begin{proposition}
The expressions  $(\divv\chih)_A=\nabb^B \chih_{AB}$, $\divv
\eta=\nabb^B\eta_B$ and $(\curll \eta)_{AB}=\nabb_A\eta_B-\nabb_B\eta_A$
verify the following equations:
\begin{align}
&(\divv \chih)_A + \chih_{AB}k_{BN}=\half (\nabb_A \trch + k_{AN} \trch) - 
\rr_{B{4}AB},
\label{Codaz}\\
&\divv\,\eta = \half\bigg(\mu +2n^{-1} N(n) \trch   -2|\eta|^2 -|\chih|^2 -2k_{AB}\chi_{AB}\bigg)
 - \half \rr_{B43A},\label{diveta}\\
&\curll\,\eta = \half\in^{AB} k_{AC} \chih_{CB} -
 \half \in^{AB}\rr_{B{43}A}.
\label{curleta}   
\end{align}
We also have the Gauss equation,
\be{gauss}
2K=\chih_{AB}\chibh_{AB}-\half \trch \trchb +R_{ABAB}
\end{equation}
\label{Nab}
\end{proposition}
We add two useful commutation formulas.
\begin{lemma}
Let $\Pi_{\und{A}}$ be an m-covariant tensor tangent to the surfaces $S_{t,u}$.
Then,
\bea
\nabb_B \ddd_{4} \Pi_{\und{A}} - \ddd_{4}\nabb_B \Pi_{\und{A}} &=&
\chi_{BC} \nabb_C \Pi_{\und{A}} - n^{-1} \nabb_B n \ddd_{4} \Pi_{\und{A}}
\label{dbpi}\\ &+& \sum_i (\chi_{A_i B} \bk_{CN} -
\chi_{BC}\bk_{A_i N} + 
\rr_{CA_i{4}B}) \Pi_{A_1..\Check{C}..A_m}.\nn
\eea
Also, for a scalar function $f$,
\begin{equation}
\label{compaper1}
\nabb_N \nabb_A f - \nabb_A\nabb_N f = -\frac 32 k_{AN}\dd_4 f - (\eta_A + k_{AN}) \dd_3 f -
(\chi_{AB}-\chib_{AB})\nabb_B f 
\end{equation}
\label{2Comm}
\end{lemma}
\begin{proof}
For simplicity we  only provide the proof of the identity \eqref{compaper1}.  The
derivation of \eqref{dbpi}  is only slightly more involved (see \cite{Ch-Kl},
\cite{Kl-Ro}). We have
$$
\nabb_N \nabb_A f - \nabb_A\nabb_N f = [N, e_A] f - (\nabb_N e_A) f=
(\dd_N e_A - \nabb_N e_A) f - (\dd_A N) f
$$
Now using the identity $N=\frac 12(e_4 -e_3)$ and the Ricci equations 
\eqref{ricciform} we can easily infer \eqref{compaper1}.
\end{proof}

\section{Special  structure of the curvature tensor $\rr$}
\label{se:specialstr}
 In this section we describe some  remarkable
decompositions\footnote{The results of this section apply to an arbitrary
 Lorentz metric $H$.} of the curvature
tensor of the metric $H$. We consider given
a system of  coordinates\footnote{This applies to the original wave coordinates
$x^\a$.}  $x^\a$ 
relative to  which
$H$ is a non degenerate
 Lorentz metric  with bounded components $H_{\a\b}$.
We define the coordinate dependent   norm 
\be{Hnorm}|\pr H|=\max_{\a,\b, \ga}|\pr_\ga H_{\a\b}|
\end{equation}
We say that  a frame $e_a, e_b, e_c, e_d$ is bounded,
with respect to our given coordinate system, if all components of 
$e_a=e_a^\a\pr_\a$  are bounded.

Consider an arbitrary bounded  frame $e_a, e_b, e_c, e_d$ and
$\rr_{abcd}$ the components of the curvature tensor relative to it.
Relative to any   system of  coordinates  we can write
\be{riemcoord}
\rr_{abcd}=e_a^\a e_b^\b e_c^\ga e_d^\de( \pr^2_{\a\ga} H_{\b\de}
+\pr^2_{\b\de} H_{\a\ga}- \pr^2_{\b\ga} H_{\a\de}-\pr^2_{\a\de} H_{\b\ga})
\end{equation}
Using our given coordinates $x^\a$ we introduce the 
flat Minkowski  metric $m_{\a\b}=\mbox{diag}(-1,1,1,1)$. We denote by
 $\nabc$ the corresponding flat  connection. Using $\nabc$ we define the following
tensor:
$$\pi(X,Y,Z)=\nabc_Z H(X, Y)$$
Thus in our local coordinates $x^\a$ we have $\pi_{\a\b\ga}=\pr_\ga H_{\a\b}$.
 \begin{proposition} Relative to an arbitrary bounded  frame $e_a, e_b, e_c, e_d$ 
we have the following decomposition:
\begin{equation}
\rr_{abcd}= \dd_a\pi_{bdc}+\dd_b\pi_{acd}-\dd_a\pi_{bcd}-\dd_b\pi_{dac} +E_{abcd}
\end{equation}
where the components of the tensor  $E$  are bounded  pointwise  by 
the square of the  first derivatives of  $ H$.
More precisely, denoting $|E|=\max_{a,b,c,d}|E_{abcd}|\approx
\max_{\a,\b,\ga,\de}|E_{\a\b\ga\de}|$, we have
\be{curv1}
|E|\les |\pr H|^2
\end{equation}
\label{Curd}
\end{proposition}
\begin{remark} It will be clear from the proof below that
we can interchange the indices $a,c$ and $b,d$ in the formula above
and obtain similar decompositions.
\end{remark}
We show that each term appearing in \eqref{riemcoord} can be expressed
in terms of a corresponding  derivative of $\pi$ plus terms of type $E$.

Consider the term $ R_1=e_a^\a e_b^\b e_c^\ga e_d^\de\pr^2_{\a\de} H_{\b\ga}$.
We show that it can be expressed in the form 
$\dd_a\pi_{bcd}$ plus terms of type $E$.  Indeed,
\beaa
\dd_a\pi_{bcd}&=&e_a(\pi_{bcd})-\pi_{\dd_ab cd}-\pi_{b\dd_a cd}-
\pi_{b c\dd_ad}\\
&=& e_a^\a \pr_\a( e_d^\de e_b^\b e_c^\ga \pr_\de H_{\b\ga})-\pi_{\dd_ab
cd}-\pi_{b\dd_a cd}-
\pi_{b c\dd_ad}\\
&=& R_1+e_a^\a \pr_\a( e_d^\de e_b^\b e_c^\ga )\pr_\de H_{\b\ga}-
\pi_{\dd_ab
cd}-\pi_{b\dd_a cd}-
\pi_{b c\dd_ad}\\
&=&R_1+ e_d^\de e_a^\a \pr_\a (e_b^\b) e_c^\ga \pr_\de H_{\b\ga}-
\pi_{\dd_ab
cd}-...
\eeaa
Now,
$$
\pi_{\dd_ab cd}=\nabc_d H(\dd_a e_b, e_c)= e_d^\de (\dd_ae_b)^\b
e_c^\ga\pr_\de H_{\b\ga}
$$
Thus,
$$
\dd_a\pi_{bcd}=R_1+e_d^\de e_c^\ga \pr_\de H_{\b\ga}\big( e_a^\a \pr_\a (e_b^\b) -
(\dd_a e_b)^\b\big)
$$
On the other hand 
\beaa (\dd_a e_b)^\b&=&<\dd_a e_b, \pr_\mu>H^{\b\mu}\\
&=&
e_a^\a\pr_\a(e_b^\b)-<e_b, \dd_a\pr_\mu>H^{\b\mu}-<e_b,\pr_\mu>e_a^\a\pr_\a
(H^{\b\mu})
\eeaa
Henceforth,
we infer that,
$$
R^{(1)}_{abcd}=\dd_a\pi_{bcd} +E^{(1)}_{abcd}
$$
with 
$$ E^{(1)}=e_d^\de e_c^\ga \pr_\de H_{\b\ga}\big(<e_b, \dd_a\pr_\mu>H^{\b\mu}+
<e_b,\pr_\mu>e_a^\a\pr_\a
(H^{\b\mu}) \big).
$$
 Since   $\dd_a\pr_\mu$ can be expressed in terms
of  the first derivatives\footnote{recall that $\dd_{\b}\pr_\mu=\Ga^\ga_{\b\mu} \pr_\ga $
with $\Ga$ the standard Christoffel symbols of $H$.} of
$H$ we conclude that
$|E^{(1)}|\les |\pr H|^2$ as desired. The other terms
in the formula \eqref{riemcoord} can be handled in precisely the same
way.
\begin{remark}  
 We will  apply proposition \ref{Curd} to
our metric $H$, wave coordinates  $x^\a$ and our canonical 
null frames. We remark that our wave coordinates are
non degenerate  relative to $H$, see \eqref{4.1}, and any canonical
null frame $e_4=(T+N), e_3=(T-N)$, $e_A$ is bounded relative
to $x^\a$.
\end{remark}
\begin{corollary}
Relative to an arbitrary frame $e_A$ on $S_{t,u}$ we have,
\be{ABCD}
\rr_{ABCD}= \nabb_A\pi_{BDC}+\nabb_B\pi_{ACD}-\nabb_A\pi_{BCD}-\nabb_B\pi_{DAC}
+E_{ABCD}
\end{equation}
with  $E$ is an error term  of the type,
$$ |E|\lesssim (|\pr H|^2+|\chi||\pr H|).$$
and, 
$$ |\pi|\les |\pr H|.$$
\label{CorABCD}
\end{corollary}
\begin{corollary}
There exists a scalar $\pi$, an S-tangent 2-tensor $\pi_{AB}$ 
and 1-form  $E_A$ such that,
the component $R_{B4AB}$ admits the decomposition
$$R_{B4AB}=\nabb_A\pi +\nabb^B\pi_{AB}+E_A.
$$
Moreover,
\beaa
|\pi|&\les& |\pr H|\\
|E|&\les &(|\pr H|^2+|\chi||\pr H|).
\eeaa

\label{corB4AB}
\end{corollary}
\begin{corollary} There exists an S-tangent vector $\pi_A$ and scalar $E$
such that
$$\in^{AB} \rr_{AB34}=\curll \pi+E$$
 and,
\beaa
|\pi|&\les& |\pr H|\\
|E|&\les &(|\pr H|^2+|\chi||\pr H|).
\eeaa
\label{corAB34}
\end{corollary}
\begin{corollary} There exist S-tangent vectors $\pi^{(1)}_A, \pi^{(2)}_A$ and scalars 
$E^{(1)},  E^{(2)}$
such that
$$
\aligned
&\de^{AB} \rr_{A43B}=\divv \pi^{(1)}+\rr + \rr_{34} + E^{(1)},\\
&\in^{AB} \rr_{A43B}=\curll \pi^{(2)} + E^{(2)}, 
\endaligned
$$
where $\rr$ is the scalar curvature.
Moreover,
\beaa
|\pi^{(1,2)}|&\les& |\pr H|\\
|E^{(1,2)}|&\les &(|\pr H|^2+|\chi||\pr H|).
\eeaa 
\label{corA4B3}
\end{corollary}
\begin{proof}
Observe that 
$
\rr_{AB}=H^{\mu\nu}\rr_{A\mu B\nu}=-\half \rr_{A3B4}-\half \rr_{A4B3}
-\de^{CD}\rr_{ACBD}.
$
Hence, since $\rr_{A3B4}=\rr_{B4A3}$, we have 
$\de^{AB}\rr_{AB}=-\de^{AB}\rr_{A4B3}-\de^{AB}\de^{CD}\rr_{ACBD},$
and therefore,
$$
\aligned
\de^{AB}\rr_{A43B}&=\de^{AB}\rr_{AB}+\de^{AB}\de^{CD}\rr_{ACBD}\\ &=
\rr+ \rr_{34}+ \de^{AB}\de^{CD}\rr_{ACBD}.
\endaligned
$$
We now appeal to corollary \ref{CorABCD} and express $\de^{AB}\rr_{A43B}$
in the form 
$$\de^{AB}\rr_{A43B}=\divv \pi^{(1)} +\rr + \rr_{34}+E^{(1)},$$
where 
\beaa
|\pi^{(1)}|&\les& |\pr H|\\
|E^{(1)}|&\les &(|\pr H|^2+|\chi||\pr H|).
\eeaa

On the other hand since $\rr_{A3B4}+\rr_{AB43}+\rr_{A43B}=0$,
we infer that $\rr_{A3B4}-\rr_{A4B3}=-\rr_{AB43}$. Thus,
$$2\in^{AB}\rr_{A43B}= -\in^{AB}\rr_{AB43}.$$
In view of corollary \ref{corAB34} we can therefore express
$\in^{AB}\rr_{A43B}$ in the form $\curll\pi^{(2)} + E^{(2)}$.
\end{proof}

\section{Strategy  of the proof of the Asymptotic  Theorem}
\label{se:strategy}
In this section we describe the main ideas in the proof of the 
Asymptotic theorem. 
\begin{enumerate}
\item {\sl Section 6}

\noindent  We start by making some primitive assumptions, which we refer to as
 \begin{itemize} 
\item Bootstrap assumptions.
\end{itemize}
 They concern 
the geometric  properties of the $C_u$  and $\stu$ foliations.
Based on this assumptions we derive further important properties,
such as
\begin{itemize}
 \item Sharp comparisons between the functions $u, r$ and $s$.
\item  Isoperimetric and Sobolev inequalities on $S_{t,u}$.
\item  Trace inequality; restriction of functions
in $H^2(\Si_t )$  to $\stu$.
\item Transport Lemma
\item Elliptic estimates on Hodge systems.
\end{itemize}
\item
 {\sl Section 7}
\noindent We recall the \textit{ background estimates} on $H=H_{(\la)}$
proved  in \cite{Einst1}. We establish further estimates of $H$ related to the surfaces 
$\stu$ and null hypersurfaces $C_u$.
\begin{itemize}
\item $L^q(\stu)$ estimates for $\pr H$ and $\ric(H)$.
\item  Energy estimates on  $C_u$.
\item  Statement of the estimate for the derivatives of $\ric_{44}(H)$.
\end{itemize}
\item
 {\sl Section 8}

\noindent Using the bootstrap assumptions 
and the results of sections 6 and 7 we provide
a detailed proof of the Asymptotics theorem.
\end{enumerate}
\section{Bootstrap assumptions and Basic Consequences}
Throughout this section we shall use only  the following background
property, see proposition 2.4 in \cite{Einst1}, of the metric $H$ in $[0,t_*]\times \rrrr^3$:
\be{someHass}
\|\pr H\|_{L_t^2L_x^\infty}\les \la^{-\f12-4\ep_0}
\end{equation}
By H\"older inequality  we also have,
\be{someHass1}
\|\pr H\|_{L_t^1L_x^\infty}\les \la^{-8\ep_0}
\end{equation}
The maximal time $t_*$ verifies the estimate $t_*\le \la^{1-8\eps_0}$.

 \subsection{Bootstrap assumptions} We start by 
constructing  the outgoing null geodesics originating from
the axis $\Ga_t$, $t\in [0, t_*]$. The geodesics emanating
from the same points $\in \Ga_t$ form the null cones $C_u$.  We define 
$\Omega^*\subset [0,t_*]\times \rrrr^3 $ to be the largest set  
properly foliated by the null cones $C_u$ with the following  properties:

{\bf A1}) \,\,\, Any point in $\Omega^*$ lies on a unique outgoing null geodesic segment
initiated from $\Ga_t$ and contained  in $\Omega^*$.

\medn {\bf A2})\,\,\, Along any fixed  $C_u$,  $ \frac{r}{s}\rightarrow 1$  as $s\rightarrow
0$.  Here  $s$ denotes the affine parameter along $C_u$, i.e. $L(s)=1$ and
$s|_{\Ga_t}=0$.  Recall  also that $r=r(t,u)$ denotes the 
radius  of $S_{t,u}=C_u\cap \Si_t$.

Moreover, the following 
bootstrap assumptions are satisfied for some $q>2$, sufficiently close to $2$ :

\medn {\bf B1})\,\,\,\, $\|\trch -\frac 2r\|_{L^2_t L^\infty_x}
\les \la^{-\frac 12-2\eps_0}$, \qquad 
$\|\chih\|_{L_t^2L_x^\infty}\les \la^{-\frac 12-2\eps_0},$
\qquad $\|\eta\|_{L_t^2L_x^\infty}\les \la^{-\frac 12 -2\eps_0},$

\medn {\bf B2})\,\,\,\,  $\|\trch - \frac 2r\|_{L^q(\stu)}\les \la^{-2\eps_0}$,\qquad
$\|\chih\|_{L^q(\stu)}\les \la^{-2\eps_0}$,
\qquad $\|\eta\|_{L^q(\stu)}\les \la^{-2\eps_0}$.

\medn

\begin{remark}
It is straightforward to check that ${\bf B1)}$ and ${\bf B2)}$ are verified
in a small neighborhood of the time axis $\Ga_t$. Indeed for each
fixed $\la$ our metrics $H_\la$ are smooth and therefore we can find 
as sufficiently small neighborhood, whose size possibly depends on $\la$,
where the assumptions  ${\bf B1)}$ and ${\bf B2)}$ hold.
\end{remark}

\begin{remark} We shall often have to estimate functions $f$ in $\Om_*$
which verify equations of the form $\frac{df}{ds}=F$ with $f=f_0$
on the axis $\Ga_t$.  According to {\bf A1})  we can express the value
of $f$ at every point $P\in \Om_*$ by the formula,
$$f(P)= f_0(P_0)+\int_{\ga}F$$
with $\ga$ the unique null geodesic in $\Om_*$ connecting the point 
$P$ with the time axis $\Ga_t$ and $P_0=\ga\cap\Ga_t$. For convenience we shall
rewrite this formula, relative to the affine parameter $s$ in the form
$$f(s)=f(0)+\int_0^s F(s') ds'.$$
It will be clear from the  context that the integral with respect to $s'$
denotes the integral along a corresponding  null geodesic $\ga$.
\end{remark}

\subsection{Comparison results}
We start with some simple  comparison\footnote{In \cite{Kl-Ro} we had
in fact $n=1$ and  $s=t-u$. In our context this is no longer
true due to the non triviality of the lapse function $n$.}  between
the affine parameter $s$ and $n(t-u)$.
\begin{lemma}
In the region $\Om_*$
$$
s\approx (t-u),\quad \mbox{i.e.},\quad s\les (t-u)\,\,\,\mbox{and}\,\,\,\,
(t-u)\les s
$$
\label{Comst}
\end{lemma}
\begin{proof}
Observe
that 
 $\frac {dt}{ds}= L(t) = T(t) = n^{-1}$
and, since $u|_{\Ga_t}=t$,
\be{1stu}
t-u = \il_{\ga} n^{-1}=\int_0^sn^{-1}(s')ds'
\end{equation}

Thus, since $n$ is bounded uniformly from 
below and above, we infer  that
 $s$ and $t-u$ are
comparable,  i.e. $s\approx t-u$. 
In particular $ s\le \la^{1-4\eps_0}$ everywhere in $\Om_*$.

\end{proof}
\begin{remark}
\label{integralst}
 The formula $\frac{ds}{dt}=n$ along $\ga$ together with
the uniform  boundedness of $n$,  used  in lemma \ref{Comst} above, 
  allows us to
estimate integrals along the  null geodesics $\ga$ as follows:
$$
|\int_{\ga} F |= |\int_0^s F(s') ds'| =|\int_0^sF(t(s'), x(s'))  ds'|=|\int_0^t(nF)(t',
x(s'(t'))dt' |\les  \|F\|_{L_t^1L_x^\infty}.
$$
 We shall make
a frequent use of this remark.
\end{remark}

In what follows we shall refine the comparison between $s$ and 
$t-u$.
\begin{lemma}
In the region $\Omega_*$
$$
n(t-u) = s\bigg(1+O(\la^{-4\eps_0}\bigg).
$$
\label{Sntu}
\end{lemma}
\begin{proof}
Consider $U= \big(n(t-u)-s\big)$ and proceed as in lemma above
by noticing that
 $\frac{du}{ds}=0$.
Therefore,
\beaa
\frac{d}{ds} U= \frac{d}{ds} \bigg(n(t-u)-s\bigg)&=&n^{-1} L (n) n(t-u)\\
&=&n^{-1} L (n)s+n^{-1} L (n)
\bigg(n(t-u)-s\bigg)
\eeaa
 Integrating  from the axis $\Ga_t$ we find,
\begin{equation}
U(s)= \int_{\ga}s' n^{-1}L(n) ds'+\int_{\ga}U(s') n^{-1}L(n) ds'
\label{stnu3}
\end{equation}
where $\ga $ is the null geodesic initiating 
on the axis $\Ga_t$ and passing through a point  $P_0$ corresponding to the value $s$. 
By Gronwall we find,
$$ U(s)\les \int_0^s s' |n^{-1}L(n)| ds'\exp{\int_0^s  |n^{-1}L(n)| ds'}.$$

According to the Remark \ref{integralst}, 
$\int_0^s n^{-1} |L(n)|\les \|\pr H\|_{L^1_t L^\infty_x}$. We can
 now make use of the inequality \eqref{someHass1}  and infer that
$$
n(t-u)= s\bigg(1+O(\la^{-8\eps_0})\bigg).
$$

\end{proof}

\begin{lemma}
The lapse function $b$ satisfies the estimate
\be{2lap}
|b(s)-n(s)|\les \la^{-8\eps_0}
\end{equation}
throughout the region $\Omega_*$.
\label{2Lap}
\end{lemma}
\begin{proof}
Integrating the transport equation \eqref{D4a}, $L(b) = - b\, \bk_{NN}$,
along the null geodesic $\ga(s)$, we infer
that,
$$
b(s)= b(0) \exp {(-\il_0^s \bk_{NN})}.
$$ 
Since $|\bk_{NN}|\les |\pr H|$, the condition \eqref{someHass1}
gives $\il_0^s |\bk_{NN}|\les \la^{-8\eps_0}$.
According to our definition $b^{-1} = T(u)$ and $u|_{\Ga_t} = t$.
Thus  $b^{-1}(0)= T(t)= n^{-1}(0)$ and therefore,
$|b(s)-n(0)|\les \la^{-8\ep_0}$. To finish the proof it only remains
to observe that 
$|n(s)-n(0)|\le \int_\ga|L(n)|\les \la^{-8\ep_0}$.
\end{proof}

Recall that the Hardy-Littlewood  maximal function\footnote{restricted to the 
interval $[0,t_*]$}  ${\cal M}(f)(t)$ of
$f(t)$ is defined by
$$
{\cal M}(f)(t) = \sup_{t_0}\frac 1{|t-t_0|} \il_{t_0}^t f(\tau)\,d\tau,
$$ 
and that,
$$\|{\cal M}(f)\|_{L_t^p}\les \|f\|_{L_t^p}$$
for any $1<p<\infty$.

\begin{lemma}
Let $a$ be a solution of the transport equation 
$$
L(a) = F
$$
Then for any point $P\in \Omega_*\cap\Si_t\cap \ga$,  where $\ga$ is the null geodesic 
initiating on the axis $\Ga_t$ at the point $P_0\in \Si_{t_0}$   and terminating at the point
$P$, we  
have the estimate
\begin{equation}
|a(P)- a(P_0)|\les s {\cal M} (\|F\|_{L^\infty_x})(t)
\label{maxim1}
\end{equation}
with $s$ the value of the affine parameter of $\ga$ corresponding to $P$.
\label{Maxim}
\end{lemma}
\begin{proof}
Integrating the equation $L(a) =\frac {da}{ds} =F$ along $\ga$ we obtain
$$
|a(P)-a(P_0)| = |\il_{\ga} F | \les \il_{t_0}^t \|F\|
_{L^\infty_x(\Si_\tau)}\,d\tau  \les (t-t_0)   {\cal M} (\|F\|_{L^\infty_x})(t)
$$
It remains to observe that $t-t_0 = t-u$  and that according 
to lemma \ref{Comst},  
$
|t-u| \les s
$
\end{proof}
Using lemma \ref{Maxim} we can now refine the conclusions
of lemmas \ref{2Lap}, \ref{Sntu}.
\begin{corollary}
\begin{align}
&b = n + s\, O\big(\MH (t)\big),\label{3lap1}\\
& n(t-u) = s + s^2 O\big ({\MH})(t)\big ),\label{sntu3}\\
&|\frac 1{n(t-u)} - \frac 1s|\les \MH(t),\label{sntu4}\\
&\|\frac 1{n(t-u)} - \frac 1s\|_{L^2_t L^\infty_x}\les \la^{-\f12 -4\eps_0}
\label{stnu5}
\end{align}
where $\MH(t)$ is the maximal function of $\|\pr H(t)\|_{L^\infty_x}$.
\label{CoSn}
\end{corollary}
\begin{proof}
The proof of \eqref{3lap1} is straightforward since 
$L(b-n)=- b \bk_{NN}-L(n)$. Now  observe that  the right hand-side 
$|b \bk_{NN}+L(n)||\les |\pr H|$ and $(b-n)|_{\Ga_t}= 0$.

Since, according to lemma \ref{Sntu},  $n(t-u)\le 2s$, the equation $L(n(t-u)-s) = n^{-1}
L(n) n(t-u)$   can be written in the form
$$
|\frac d{ds}\big (n(t-u) -s\big )| \les s |\pr H|
$$
Thus with the help of lemma \ref{Maxim} we obtain
$$
|n(t-u) -s|\les s^2 \MH
$$
The inequality \eqref{sntu4} is an immediate consequence of \eqref{sntu3} and
lemma \ref{Sntu}.
The estimate  \eqref{stnu5} follows from \eqref{sntu4}, 
\eqref{someHass}, and the $L^2$ estimate for the Hardy-Littlewood
maximal function.
\end{proof}
   
We shall now compare the values of the parameters $s$ and 
$r= \frac 1{4\pi} A^{\half} (\stu)$ at a point $P\in \stu$.
\begin{lemma} The identity 
$$
r= s\bigg (1+O(\la^{-6\eps_0})\bigg),
$$
holds throughout the region $\Omega_*$. In particular this implies that
$$
2\pi s^2\le A(t,u)\le 8\pi s^2
$$
with $A(t,u)$ the area of $S_{t,u}$.
\label{Rs}
\end{lemma}
\begin{proof}
Similarly to \eqref{randt}, we have
$$
L(r) = \frac r{2} \overline {\trch} = \frac 1{8\pi r}\il_{\stu}\trch
$$
Using the identity $A(\stu) = 4\pi r^2$, we obtain
\be{eqrs}
\frac {dr}{ds} = 1 +  \frac 1{8\pi r}\il_{\stu}\bigg (\trch - \frac 2r\bigg )
\end{equation}
Integrating along the null geodesic $\ga$ passing through the
point $P=P(s)$ \footnote{Observe that according to ${\bf A2)}$,
$(r-s)\to 0$ as $s\to 0$ along $C_u$}  we have
\begin{equation}
\begin{split}
|r(P) - s|&\les \il_{\ga} \frac 1{r} \il_{\stu}\bigg (\trch - \frac 2r\bigg )\le
4\pi\int_\ga  r \|\trch - \frac 2r\|_{L_x^\infty}\\
&\les \int_\ga  (r-s') \|\trch - \frac 2r\|_{L_x^\infty}+\int_\ga  s' \|\trch - \frac
2r\|_{L_x^\infty}
\end{split}
\label{crs}
\end{equation}
Thus by Gronwall,  and the bootstrap  estimate ${\bf B1}$), 
$$\|tr\chi-\frac 2r\|_{L_t^1L_x^\infty}\les \la^{\f12-4\ep_0}\|tr\chi-\frac
2r\|_{L_t^2L_x^\infty}\les \la^{-6\ep_0}$$
 we infer that,
$|r-s|\les s\la^{-6\ep_0}.$

\end{proof}
Having established that $r\approx s$ we shall  now derive more  refined comparison
estimates involving $\trch-\frac{2}{s}$  and its iterated  maximal functions. 
These will be needed later on in section \ref{se:trest} where $\trch-\frac{2}{s}$ rather
than $\trch-\frac 2r$ appears naturally.
\begin{corollary}
\begin{align}
&|r-s|\les s^2 {\cal M}^3 (\|\trch -\frac 2s\|_{L^\infty_x}),\\
&|r-s|\les s^{\frac 32} \|\trch - \frac 2s\|_{L^2_t L^\infty_x},
\end{align}
Here, ${\cal M}^k$ is the $k$-th maximal function. 
Moreover,
\begin{align}
&|\trch -\frac 2r|\les |\trch -\frac 2s| + {\cal M}^3 
(\|\trch -\frac 2s\|_{L^\infty_x}),\label{2cors1}\\
&\|\trch -\frac 2r\|_{L^2_t L^\infty_x} \les   
\|\trch -\frac 2s\|_{L^2_t L^\infty_x}\label{cors1},\\
&\|\trch -\frac 2r\|_{L^q(\stu)} \les   (1+ r^{\frac 2q-\frac 12})
\|\trch -\frac 2s\|_{L^2_t L^\infty_x},\label{cors2}\\
&\|\frac 2r -\frac 2{n(t-u)}\|_{L^2_t L^\infty_x} \les   
\|\trch -\frac 2s\|_{L^2_t L^\infty_x} + \la^{-\f12 -4\eps_0} \label{cors3}
\end{align}
\label{Comp}
\end{corollary}
\begin{proof}
We write the transport equation for $r$ in the following form:
\begin{equation}
L(r) = \frac 1{8\pi r} \il_{\stu} (\trch - \frac 2s) + \frac 1{8\pi r} \il_{\stu}\frac 2s
\label{comp1}
\end{equation}
Differentiating $\,\il_{\stu}\frac 2s$ we obtain
\be{l2s}
L\bigg(\il_{\stu}\frac 2s\bigg ) = \il_{\stu} (\frac 2s\trch -\frac 2{s^2}) =
 \il_{\stu}\frac 2s(\trch - \frac 2s) +  \il_{\stu}\frac 2{s^2}
\end{equation}
Furthermore,
$$
L\bigg(\il_{\stu}\frac 2{s^2}\bigg ) =
2 \il_{\stu}\frac 1{s^2}(\trch - \frac 2s)
$$
Since $s-r\to 0$ as $r\to 0$, we have $\il_{\stu}\frac 2{s^2}\to 8\pi$.
Using lemmas \ref{Rs} and \ref{Maxim} we infer that
$$
\il_{\stu}\frac 2{s^2} = 8 \pi + s {\cal M}\big(\|\trch -\frac 2s\|_{L^\infty_x}\big)
$$
Integrating \eqref{l2s} and using lemma \ref{Maxim} once more we obtain
$$
\il_{\stu}\frac 2s = 8\pi s +
s^2 {\cal M}^2\big(\|\trch -\frac 2s\|_{L^\infty_x}\big) + 
s^2 {\cal M}\big(\|\trch -\frac 2s\|_{L^\infty_x}\big) 
$$
Again, according to lemma \ref{Rs},  $r\approx s$. Thus returning to 
\eqref{comp1} 
$$
L(r) = \frac sr +  \frac 1{8\pi r} \il_{\stu} (\trch - \frac 2s) + 
s {\cal M}^2\big(\|\trch -\frac 2s\|_{L^\infty_x}\big) 
$$
or, equivalently,
$$
L(r^2) = 2 s + \frac 1{4\pi } \il_{\stu} (\trch - \frac 2s) + rs
{\cal M}^2\big(\|\trch -\frac 2s\|_{L^\infty_x}\big) 
$$
Integrating with the help of lemma \ref{Maxim} we infer that,
$$
r^2 = s^2 + s^3{\cal M}^3\big(\|\trch -\frac 2s\|_{L^\infty_x}\big)
+ s^3  {\cal M}\big(\|\trch -\frac 2s\|_{L^\infty_x}\big) 
$$
It then follows that
\be{rsM}
r=s +  s^2{\cal M}^3\big(\|\trch -\frac 2s\|_{L^\infty_x}\big)
\end{equation}
Observe that if during each integration along $\ga$ we 
used H\"older inequality instead of the bounds involving
maximal functions, we would have the estimate
\be{rs2t}
r=s +  s^{\frac 32}\|\trch -\frac 2s\|_{L^2_t L^\infty_x}
\end{equation}
This estimate can be used effectively to compare $r$ and $s$ 
on a single surface $\stu$ while \eqref{rsM} works well with the
norms involving integration in time.
Thus, we infer from from \eqref{rsM} that
\begin{align}
&|\frac 2r -\frac 2s|\les {\cal M}^3\big(\|\trch -\frac 2s\|_{L^\infty_x}\big),\label{10cors}\\
& \|\frac 2r -\frac 2s\|_{L^2_t L^\infty_x} \les
\|\trch -\frac 2s\|_{L^2_t L^\infty_x}\label{cors4}
\end{align}
In addition, \eqref{rs2t} implies that
\be{rsq}
\|\frac 2r -\frac 2s\|_{L^q(\stu)} \les r^{\frac 2q-\frac 12} 
\|\trch -\frac 2s\|_{L^2_t L^\infty_x}
\end{equation}
Inequalities \eqref{2cors1}-\eqref{cors2} follow from the identity 
$\trch-\frac 2r = \trch -\frac 2s  +\frac 2r -\frac 2s$ and 
\eqref{10cors}-\eqref{rsq}.
Finally, \eqref{cors3} follows from \eqref{10cors} and  \eqref{stnu5}.
\end{proof}

\begin{remark}

Observe that the equation \eqref{eqrs} and lemma \ref{Maxim}
also give the estimate
$$
|r-s|\les s^2 {\cal M}\big(\|\trch - \frac 2r\|_{L^\infty_x}\big)(t).
$$
Thus with the help of the bootstrap assumption ${\bf B1)}$ and the $L^2$
estimate for the maximal function we infer that,
\begin{equation}
\begin{split}
\|\trch - \frac 2s\|_{L^2_t L^\infty_x} &\les 
\|\trch - \frac 2r\|_{L^2_t L^\infty_x} + \|\frac 2r -\frac 2s\|_{L^2_t L^\infty_x}\\
&\les 2  \|\trch - \frac 2r\|_{L^2_t L^\infty_x}\les \la^{-\frac 12 -2\eps_0}
\end{split}
\label{rs2}
\end{equation}
Moreover, since $r\approx s$, the equation \eqref{eqrs}, H\"older 
inequality and the bootstrap assumption ${\bf B2)}$ also
imply that
$$
|r-s|\les \int_{\ga} r^{1-\frac 2q} \|\trch -\frac 2r\|_{L^q(\stu)}
\les \la^{-2\eps_0} s \, r^{1-\frac 2q}
$$
Using the bootstrap assumption ${\bf B2)}$ once again we infer that
\begin{equation}
\begin{split}
\|\trch -\frac 2s\|_{L^q(\stu)} &\les \|\trch -\frac 2r\|_{L^q(\stu)}
+ \|\frac 2r -\frac 2s\|_{L^q(\stu)}\\ &\les \la^{-2\eps_0} + 
\la^{-2\eps_0}  \|r^{-\frac 2q}\|_{L^q(\stu)}\les \la^{-2\eps_0}
\end{split}
\label{rs1}
\end{equation}
Estimates \eqref{rs1}, \eqref{rs2} indicate that the bootstrap assumptions
${\bf B1)}$, ${\bf B2)}$ also hold for $(\trch -\frac 2s)$.
\label{RRs}
\end{remark}

\subsection{Isoperimetric, Sobolev inequalities and   transport lemma}

We consider now the foliation induced by $S_{t,u}$ on $\Si_t\cap\Omega_*$.
Relative to this foliation the induced metric $h$  on $\Si_t$ takes the form
$$ h=b^2 du^2 +\ga_{AB} d\phi^A d\phi^B$$
where $\phi^A$ are local coordinates on $S^2$.
We state below a proposition concerning the trace and isoperimetric
inequalities on $\Si_t\cap \Om_*$. The proposition requires 
a very weak assumption on the metric $h$,
in fact we only need
\be{assisop}
\big (\sup_{\Omega_*} r^{\frac{1}{2}\ep}\big)
\|\nab^{\frac{3}{2}+\epsilon} h\|_{L^2(\Si_t)}\le \La_0^{-1}
\end{equation}
for some large constant $\La_0>0$ and an arbitrarily small $\eps>0$.
In this and the following subsection we shall assume a slightly stronger 
property that
\be{assisop'}
\big (\sup_{\Omega_*} r^{\frac{1}{2}\ep}\big)
\|\nab^{\frac{1}{2}+\epsilon} \pr H\|_{L^2(\Si_t)}\le \La_0^{-1}
\end{equation}
\begin{remark} The assumption  \eqref{assisop'} is easily satisfied by
our families of metrics $H=H_{(\la)}$, see remark \ref{re:isop}.
\end{remark}
\begin{proposition}Let $S_{t,u}$ be a fixed surface in  $\Si_t\cap \Om_*$.

\medn
{\bf i.}\,\,\,For any smooth function $f: S_{t,u}\to \Bbb R$ we have 
the following isoperimetric 
inequality:
\begin{equation}
\Bigl (\il_{S_{t,u}} |f|^2\,\Bigr )^{\half }
\les \il_{S_{t,u}} (|\nabb f| + \frac 1r |f|).     
\label{asob3}
\end{equation}
\medn
{\bf ii.}\,\,\, The following Sobolev inequality holds on $S_{t,u}$: 
for any $\de\in (0,1)$ and
$p$ from the interval $p\in (2,\infty]$

\begin{equation}
\begin{split}
\sup_{S_{t,u}} |f|\,\les&\, r^{\frac {\eps(p-2)}{2p+\de(p-2)}}
\bigl (\il_{S_{t,u}} (|\nabb f|^2 + r^{-2} |f|^2)
 \bigr)^{\frac 12-\frac{\de p}{2p+\de(p-2)}}\\
&\bigg[{\il_{S_{t,u}} (|\nabb f|^p + r^{-p} |f|^p)}
\bigg]^{\frac{2\de}{2p+\de(p-2)}}, 
\end{split}    
\label{asob4}
\end{equation}

\medn {\bf iii.}\,\,\,
Consider an arbitrary function $f:\Si_t\rightarrow \rr$ such that 
$f\in H^{\f12+\ep}(\rr^3)$. The following trace inequality holds true:
\be{traceineq}
\|f\|_{L^2(S_{t,u})}\les\|\pr^{\f12 +\ep}f\|_{L^2(\Si_t)}+\|\pr^{\f12
-\ep}f\|_{L^2(\Si_t)}.
\end{equation}
More generally, for any $q\in [2,\infty)$
\be{traceineq2}
\|f\|_{L^q(S_{t,u})}\les\|\pr^{\frac 32 -\frac 2q +\ep}f\|_{L^2(\Si_t)}+\|\pr^{\frac 
32 -\frac 2q -\ep}f\|_{L^2(\Si_t)}.
\end{equation}
Also, considering the region $\Omega_*(\frac 14 r,r)=\cup_{\frac 1 4 r\le \rho\le
r}S_{t,u(\rho)}$,
 where $r=r(t,u)$, we have the following:
\be{traceineq22}
\|f\|_{L^2(S_{t,u})}^2\le\|N(f)\|_{L^2(\Omega_*(\frac 14 r,r))}\|f\|_{L^2(\Omega_*(\frac 14
r,r))}+\frac{1}{r}
\|f\|_{L^2(\Omega_*(\frac 14 r,r))}.
\end{equation}
\label{Triso}
\end{proposition}

Finally we state below,
\begin{lemma}[Transport Lemma]
Let $\Pi_{\und {A}}$ be an S-tangent  tensorfield 
verifying the following transport equation
with $\si >0$:
$$
\ddd_{4} \Pi_{\und{A}} + \si \trch \Pi_{\und{A}} = F_{\und{A}}.
$$
Assume that the point $(t,x)=(t,s,\omega)$ belongs to the domain $\Omega_*$.
If  $\Pi$ satisfies the initial condition  $s^{2\si} \Pi_{\und{A}} (s)\to 0$ 
as $s\to 0$,
then
\begin{equation}
|\Pi(t,x)|\le 4 \|F\|_{L^1_t L^\infty_x}.
\label{tran2}
\end{equation}
In addition, if $\si\ge \frac 1q$ and $\Pi$ satisfies the initial condition

$r^{2(\si -\frac 1q)}\|\Pi\|_{L^q(S_{t,u})}\to 0$ as $r\to 0$,
then on each surface $S_{t,u}\subset \Omega_*$ 
\begin{equation}
\|\Pi\|_{L^q(S_{t,u})} \les \frac{1} {r(t)^{2(\si-\frac 1q)}} \il_u^t 
r(t') ^{2(\si-\frac 1q)} 
\|F\|_{L^q(S_{t',u})}\,dt'
\label{tran3}
\end{equation}
Finally, if $\Pi$ is a solution of the transport equation
$$
\ddd_{4} \Pi_{\und{A}} + \si \trch \Pi_{\und{A}} = \frac 1r F_{\und{A}}.
$$
verifying the initial condition $s^{2\si} \Pi_{\und{A}} (s)\to 0$ with some $\si>\frac 12$,
then 
\begin{equation}
|\Pi(t,x)|\le 4 {\mathcal M}(\|F\|_{L^\infty_x})(t).
\label{20tran}
\end{equation}
\label{Tran}
\end{lemma}
\begin{proof}
The proof of \eqref{tran2}-\eqref{tran3} is straightforward. For a similar 
version see \cite{Kl-Ro}. Estimate \eqref{20tran} can be proved in the same
manner as \eqref{maxim1} of lemma \ref{Maxim}.
\end{proof}

\subsection{Elliptic estimates}

Next we establish a proposition concerning the $L^2$ estimates
of Hodge systems on the surfaces $S_{t,u}$. They are similar to
 the estimates of lemma 5.5 in \cite{Kl-Ro}. We need however to make
an important modification based on the corollary \ref{CorABCD}.

\begin{proposition}
Let $\xi$ be an $m+1$ 
covariant, totally symmetric tensor, a 
solution of the 
Hodge system on the surface $S_{t,u}\subset\Omega_*$
$$
\aligned
&\divv \xi = F, \\
&\curll \xi = G, \\
&\tr \xi = 0.
\endaligned
$$

 Then $\xi$ 
obeys the estimate
\begin{equation}
\il_{S_{t,u}} |\nabb \xi |^2 + \frac{m+1}{2r^2}|\xi|^2 \le 2 \il_{S_{t,u}} 
\{|F|^2 + |G|^2 \}.
\label{2hodge4}
\end{equation}
\label{2Hodge}
\end{proposition}
\begin{proof}
Using the standard Hodge theory, see  theorem5.4 in \cite{Kl-Ro} or chapter 2
in \cite{Ch-Kl}, we have 
\begin{equation}
\il_{S_{t,u}} |\nabb \xi |^2 + (m+1) K |\xi|^2 = \il_{S_{t,u}} \{|F|^2 + |G|^2\}.
\label{hodge4}
\end{equation}
The Gauss curvature $K$ of the 2-surface $S_{t,u}$ can be expressed as
follows: 
$$
K=\frac 14 (\trch)^2 +\f12 \trch \tr k + \f12 \chih\cdot \chibh + \f12 R_{ABAB}
$$
Thus it follows from corollary \ref{CorABCD} that
$$
K-r^{-2} = \nabb_A \Pi_A + E
$$
where the tensor $\Pi$ and the error term $E$, relative to the standard
coordinates $x^\a$, obey the pointwise estimates $|\Pi|\les |\pr H|$ and
$|E|\les (|\pr H|^2 + |\chih|^2 + |\chi||\pr H|)$.
Then we have
\be{elr}
\il_{S_{t,u}} |\nabb \xi |^2 + \frac{m+1} {r^{2}} |\xi|^2 \le 
\il_{S_{t,u}} \{|F|^2 + |G|^2 + (m+1) (\nabb_A \Pi_A + E)|\xi|^2\}.
\end{equation}
Integrating the term $\il_{S_{t,u}} \nabb_A \Pi_A |\xi|^2$ by parts we obtain
for all sufficiently large $p$, $\frac 12=\frac 1p + \frac 1q$, 
$$
\il_{S_{t,u}} \nabb_A \Pi_A |\xi|^2 = - 2\il_{S_{t,u}} \Pi_A 
\nabb_A\xi\cdot\xi \les \|\nabb \xi\|_{L^2(S_{t,u})} \|\xi\|_{L^p(S_{t,u})}
\|\Pi\|_{L^q(S_{t,u})}.
$$
The isoperimetric inequality implies that for $2\le p <\infty $
$$
\|\xi\|_{L^p(S_{t,u})}\les r^{\frac 2p}\bigg (\|\nabb \xi\|_{L^2(S_{t,u})} +
r^{-1}\|\xi\|_{L^2(S_{t,u})}\bigg )
$$
We also deduce from the trace inequality that
$$
\|\Pi\|_{L^q(S_{t,u})}\les \|\pr H\|_{L^q(S_{t,u})}
\les \|\pr^{(\frac 32+1- \frac 2q+\eps)} H\|_{L^2(\Si_t)}+
\|\pr^{(\frac 32+1-\frac 2q-\eps)} H\|_{L^2(\Si_t)}.
$$
Thus the smallness condition 
$$ 
r^{1- \frac 2q} \|\pr^{(\frac 32+ 1- \frac 2q +\eps)} H\|_{L^2(\Si_t)}\le \La_0^{-1}
$$
ensures that we can absorb the term
$(m+1)\il_{S_{t,u}} \nabb_A \Pi_A |\xi|^2$  on the left hand-side of
\eqref{elr}. For large $p$ the above condition coincides with \eqref{assisop}.

It remains to estimate $\int_\stu E|\xi|^2$.
The most dangerous term is $\int_\stu |\chih|^2 |\xi|^2$.
Applying the H\"older inequality we infer that,
$$
\int_\stu |\chih|^2 |\xi|^2\les \|\xi\|^2_{L^p(\stu)} \|\chih\|^2_{L^q(\stu)}.
$$
Using the isoperimetric inequality once more, we conclude that
we need  a smallness condition on $r^{1-\frac{2}{q}}\|\chih\|_{L^q(\stu)}$
for some $q>2$.  This is guaranteed by  our bootstrap assumption {\bf B2}).
\end{proof}
We shall next formulate a version of the Calderon-Zygmund theorem
for the above type of Hodge systems.
\begin{proposition}
Let $\xi$ be an $2$ 
covariant, traceless  symmetric tensor, verifying  the 
Hodge system on the surface $S_{t,u}\subset\Omega_*$
$$
\divv \xi = \nabb \nu+e$$
for some scalar  $\nu$ and 1-form $e$. Then,
\be{caldzygm1}
\|\xi\|_{L^q(\stu)}\lesssim \|\nu\|_{L^q(\stu)}  +  \|e\| _{L^p(\stu)}
\end{equation}
where $\frac{1}{p}=\frac{1}{2}+\frac{1}{q}$.
\label{caldzygm}

Also\footnote{The term $\|r \nabb \nu\|_{L^\infty(\stu)}$ can be in fact replaced 
by $\|r \nabb \nu\|_{L^r(\stu)}$ for $r>2$},
\be{caldzygm2}
\|\xi\|_{L^\infty(\stu)}\les \|\nu\|_{L^\infty(\stu)}
\logp(r \|\nabb \nu\|_{L^\infty(\stu)}) +
r^{1-\frac{2}{p}}\|e\| _{L^p(\stu)}
\end{equation}
for any $p>2$, where $\logp z = \log (2+|z|)$.

Similar estimates hold in the case when $\xi$ is a 1-form verifying 
the Hodge system
$$
\aligned
&\divv \xi = \divv \nu_1 + e_1,\\
&\curll \xi = \curll \nu_2 + e_2
\endaligned
$$
for some 1-forms $\nu=(\nu_1,\nu_2)$ and 
scalars $e=(e_1,e_2)$. 
\label{CZestimates}
\end{proposition}

\section{Properties of the  metric $H$ and its curvature tensor $\rr$}
\subsection{Background estimates}
We start by recalling the \textit{background} estimates on the family
of the Lorentz metrics $H=H_{(\la)}$ proved in \cite{Einst1}, see 
proposition (??).

Metric $H$ admits the canonical decomposition
$$
H=-n^2 dt^2 + h_{ij}(dx^i + v^i dt)\otimes (dx^j+v^j dt)
$$
and satisfies the following estimates on the time interval
$[0,t_*]$ with $t_* \le\la^{1-8\eps_0}$:
\be{psith}
c|\xi|^2\le h_{ij}\xi^i\xi^j \le c^{-1}|\xi|^2,\quad  n^2-|v|^2_h\ge c>0, \quad
|n|,|v|\le c^{-1}
\end{equation}
\begin{align}
&\|\pr^{1+m}H\|_{L^1_{[0,t_*]} L_x^\infty}\les \la^{-8\eps_0},
\label{ash1}\\
&\|\pr^{1+m}H\|_{L^2_{[0,t_*]} L_x^\infty}\les
\la^{-\frac 12 -4\eps_0},
\label{ash2}\\
&\|\pr^{1+m} H\|_{L^\infty_{[0,t_*]} L_x^\infty}\les 
\la^{-\frac 12-4 \eps_0},
\label{ash3}\\
&\|\nab^{\frac 12+m}(\pr H)\|_{L^\infty_{[0,t_*]} L_x^2}\les \la^{-m}\quad
\mbox{for}\quad  -\frac 12\le m \le \frac12 +4\ep_0\label{ash5}\\
&\|\nab^{\frac 12+m}(\pr^2 H)\|_{L^\infty_{[0,t_*]} L_x^2}\les 
\la^{-\f12 -4\ep_0}\quad
\mbox{for}\quad  -\frac12 + 4\ep_0 \le m  \label{ash5'}\\
&\|\nab^m \big({H^{\a\b}}\pr_\a\pr_\b H\big)\|_{L^1_{[0,t_*]} L_x^\infty}\les 
\la^{-1-8\eps_0},
\label{ash4}\\
&\|\nab^m (\nab^{\f12}\rr_{\a\b}(H))\|_{L^2_x}\les \la^{-1},
\label{ash800}\\
&\|\nab^m \rr_{\a\b}(H)\|_{L^1_{[0,t_*]} L_x^\infty}\les 
\la^{-1- 8\eps_0 }.
\label{ash6}
\end{align} 

\begin{remark}
\label{re:isop}
The inequality \eqref{ash2} with $m=0$ is consistent with the property
\eqref{someHass}, which we have used throughout section 6. 
Moreover, since in the region $\Omega_*$ the radius $r$ of the
surfaces $\stu$ does not exceed $\la^{1-8\eps_0}$, we have, 
according to \eqref{ash5},
$$
r^{\frac 12 \eps}\|\nab^{\frac 12+\eps}(\pr H)\|_{L^\infty_{[0,t_*]} L_x^2}\les 
\la^{(\frac 12 -4\eps_0)\eps} \la^{-\eps}\le \la^{-\frac 12 \eps}.
$$
This verifies the condition \eqref{assisop'}.
\end{remark}

\subsection{$L^q(\stu)$ estimates} 
The trace inequality \eqref{traceineq2} of proposition \ref{Triso}
allows us to derive the $L^q(\stu)$ estimates on the metric $H$ from
\eqref{ash5}.  

\begin{proposition}
For any $q$ in the interval $2\le q\le 4$
\begin{equation}
\|\pr H\|_{L^q(\stu)} \les  \la^{\frac 2q- 1-8(\frac 2q-\f12)\eps_0} 
\label{hq}
\end{equation}
In addition, 
\begin{equation}
\|\ric(H)\|_{L^p(\stu)}\les \la^{\frac 2p- 2- 8(\frac 2p - 1) \eps_0}
\label{hq2}
\end{equation}
for $p\in [1,2]$.
\label{Hq}
\end{proposition}
\begin{proof}

Since $q\le 4$, by H\"older inequality
$$
\|\pr H\|_{L^q(\stu)}\les r^{\frac 2q - \f12} \|\pr H\|_{L^4(\stu)}
\les \la^{(\frac 2q-\f12)(1-8\eps_0)}   \|\pr H\|_{L^4(\stu)}
$$
Using the trace estimate \eqref{traceineq2} we infer that
$$
\|\pr H\|_{L^q(\stu)}\les \la^{(\frac 2q-\f12)(1-8\eps_0)}   
\|\pr H\|_{\dot H^1(\rr^3)}\les \la^{\frac 2q- 1-8(\frac 2q-\f12)\eps_0} 
$$
where we have used $\|\pr H\|_{\dot H^1(\rr^3)}\les \la^{-\f12}$ from
\eqref{ash5}.
The inequality \eqref{hq2} follows similarly from the trace theorem
and \eqref{ash800}.
\end{proof}

\subsection{Energy estimates on $C_u$} In this subsection we shall
derive energy estimates, along the null hypersurfaces $C_u$, for
tangential derivatives of the first derivatives of the rescaled
metric 
\be{rescG}
G(t,x)=\gg(\frac t\la, \frac x\la)
\end{equation}
Recall that  the original space time Einstein metric $\gg$, 
 verifies $\rr_{\mu\nu}(\gg)=0$. 
In addition, since  our coordinates $x^\a$ satisfy the wave 
coordinate condition  \eqref{I2}, the  metric $\gg$ satisfies the
 quasilinear wave 
equation 
\be{wavemap}\gg^{\a\b} \pr_\a\pr_\b \gg_{\mu\nu} =N_{\mu\nu} (\gg,\pr \gg).
\end{equation}
We have  also defined the truncated 
$ \gg_{<\la} =\sum_{\mu<\f12  \la} P_\mu \gg$ and, by rescaling, 
$$
H(t,x) = \gg_{<\la}(\frac {t}\la, \frac {x}\la).
$$
 our background
metric.  Similarly, for a dyadic $\mu \ge \f12$ 
we can define 
$$
G^{(\mu)}(t,x) = P_{\mu\la}\gg  (\frac {t}\la, \frac {x}\la)
$$  
Observe that $H$ has frequencies $\le 1 $ and 
$\Gmu$ is localized to the frequencies of size $\mu$
which can not fall below $\f12$.

We now formulate a basic energy estimate on the null 
cones $C_u$ for $H$ and $\Gmu$.
\begin{definition}
\label{tangentialnotation} Given a scalar function  $F$ in $\Om_*$
we denote by 
 $D_*F$   the $C_u$ tangential derivatives  of $F$. More precisely,
$D_*F=(\nabb F, L F)$. We shall use this notation for the components
of the metrics $H$ and $G$ relative to our fixed system of coordinates.
We also use this notation applied to 
  all components of the derivatives $\pr H$ and $\pr G$. Thus
$|D_*\pr H|=\sum_{\a,\b,\ga}|D_*\pr_\ga H_{\a\b}|$
\end{definition}
\begin{proposition}
The following estimates hold in the region $\Omega_*$:
\be{cuH}
\|D_*\pr H\|_{L^2(C_u)} \les \la^{-\f12},\qquad
\|D_* H\|_{L^2(C_u)} \les \la^{\f12}
\end{equation}
In addition, for the functions $\Gmu$ defined above
\begin{equation}
\begin{split}
&\|D_*\pr \Gmu\|_{L^2(C_u)} \les \mu^{\frac 12 - 4\eps_0} \la^{-\frac 12 - 4\eps_0},\\
&\|D_* \Gmu\|_{L^2(C_u)}\les \min\{\mu^{-1-4\eps_0} \la^{-\frac 12-4\eps_0},
\,\mu^{-\frac 12-4\eps_0} \la^{-1-4\eps_0}\}
\end{split}
\label{10cuh}
\end{equation}
\label{CuH}
\end{proposition}
The following result can be deduced from propositions \ref{CuH}, \ref{Curd}.
\begin{corollary}
Any component of the curvature $\rr_{abcd}=\rr(e_a,e_b,e_c,e_d)$ with  
vectorfields $e_a,e_b,e_c$ varying between $L, e_A, A=1,2$,
obeys the energy estimates on $C_u$:
$$
\|\rr_{abcd}\|_{L^2(C_u)}\les \la^{-\f12}
$$
In particular,
$$
\|\rr_*\|_{L^2(C_u)}:=\sum_{A,B,C,D}\|\rr_{ABCD}\|_{L^2(C_u)} +
\|\rr_{ABC4}\|_{L^2(C_u)}+
\|\rr_{B43A}\|_{L^2(C_u)}\les \la^{-\f12}
$$
\label{CuR}
\end{corollary}
\begin{proof}{\bf of proposition \ref{CuH}}

Metric $\gg$ is a $H^{2+\ga}$ solution
of the Einstein equation. Thus after rescaling and taking into account 
 $\ga>5\eps_0$, we infer that in addition to the estimates 
\eqref{ash1}-\eqref{ash4} for $H$, we also have 
\begin{equation}
\|\pr^{1+m} \Gmu\|_{L^\infty_t L^2_x} \les \la^{-\frac 12 -4\eps_0}\mu^{m-1-4\eps_0},
\quad{\text{for}}\quad m=0,1
\label{2cuh2}
\end{equation}

We shall make use of the rescaled version of lemma {?} in \cite{Einst1} to derive the 
equations for $H$ and $\Gmu$. 
\begin{equation}
H^{\a\b} \pr_\a\pr_\b H = F, \qquad H^{\a\b} \pr_\a\pr_\b \Gmu = F_\mu,
\label{3cuh1}
\end{equation}
with the right hand-sides $F$, $F_\mu$ obeying the estimates
\begin{alignat}{2}
&\|F\|_{L^1_t L^2_x} \les \la^{\f12},&\qquad &
\|\pr F\|_{L^1_t L^2_x} \les \la^{-\f12},\label{2cuh3},\\
&\|F_\mu\|_{L^1_t L^2_x} \les \mu^{-4\eps_0} \la^{-\frac 12-4\eps_0},
&\qquad &\|\pr F_\mu\|_{L^1_t L^2_x} \les \mu^{1-4\eps_0} \la^{-\f12-4\eps_0}
\label{2cuh4}
\end{alignat}

We shall use the generalized energy identity with the vectorfield
$T$ in the region $M_{t_0,t,u}$ bounded by the cone $C_u$ and the time slices 
$\Si_{t_0}$, $\Si_t$ intersecting $C_u$.
The vectorfield $L$ is orthogonal, in the sense of the Lorentzian
metric $H$, to the cone $C_u$. Thus 
$$
\il_{C_u} Q[H](T,L) + \il_{\Si_{t_0}}  Q[H](T,T) = 
\il_{\Si_{t_0}}  Q[H](T,T) - \il_{M_{t_0,t,u}} \bigg(Q^{\a\b}[H] ^T\pi_{\a\b}
+ F T(H)\bigg)
$$
with the energy-momentum tensor 
$$
Q[f]_{\a\b} = \pr_\a f \pr_\b f -\frac 12 H_{\a\b} 
(\pr_\nu f \pr^\nu f)
$$ 
and the deformation tensor 
$ \piT_{\a\b}= {\cal L}_T H$ of the vectorfield $T$.
A similar identity also holds for $G^\mu$.
According to \eqref{lit} and \eqref{3.2} the components of the 
deformation tensor $^T\pi $ can be described as follows:
$$
\piT_{ij} = - 2 k_{ij}, \qquad \piT_{i0}= n^{-1} \pr_i n,
\qquad \piT_{00}= 0
$$
Thus the deformation tensor  $|\piT|\les |\pr H|$, and by \eqref{ash1} 
obeys the estimate 
\begin{equation}
\|\piT\|_{L^1_t L^\infty_x}\les \la^{-4\eps_0}
\label{2cuh1}
\end{equation}
Observe that 
$$
\aligned
\ & Q[H](T,L) = \frac 12 (L H)^2 + \frac 12 |\nabb H|^2 = 
\frac 12 |D_*H|^2  ,\\
\ & Q[H](T,T) = \frac 12 (T H)^2 + \frac 12 |\nab H|^2 = 
\frac 12 |\pr H|^2
\endaligned
$$
In addition, $|Q_{\a\b}(f)| \le 2|\pr f|^2$.
Thus, using \eqref{ash5}, \eqref{2cuh3}, and \eqref{2cuh1}, we obtain
$$
\aligned
\il_{C_u}  |D_*H|^2 &\le \il_{\Si_{t_0}}  |\pr H|^2  + 
4\il_{M_{t_0,t,u}} \bigg (| \piT | \,|\pr H|^2 + |F|\,|\pr H|\bigg)\\ &\les
\|\pr H\|^2_{L^\infty_t L^2_x} + \| \piT \|_{L^1_t L^\infty_x} 
\|\pr H\|^2_{L^\infty_t L^2_x} +   \| F \|_{L^1_t L^2_x}
\|\pr H\|_{L^\infty_t L^2_x}\les \la
\endaligned
$$
Similarly,
$$
\aligned
\il_{C_u}  |D_*G^\mu|^2 &\le \il_{\Si_{t_0}}  |\pr G^\mu|^2  + 
4\il_{M_{t_0,t,u}} \bigg (| ^T\pi | \,|\pr G^\mu|^2 + |F_\mu|\,|\pr G^\mu|\bigg)\\ 
&\les \|\pr G^\mu\|^2_{L^\infty_t L^2_x} + \| ^T\pi \|_{L^1_t L^\infty_x} 
\|\pr G^\mu\|^2_{L^\infty_t L^2_x} +   \| F_\mu \|_{L^1_t L^2_x}
\|\pr G^\mu\|_{L^\infty_t L^2_x}\\ &\les \min\{\mu^{-2-8\eps_0} \la^{-1-8\eps_0},
\mu^{-1-8\eps_0} \la^{-2-8\eps_0}\}
\endaligned
$$
To get the estimates for $D_*\pr H$ and $D_*\pr G^\mu$ we differentiate 
the equations \eqref{3cuh1}. Commuting the derivative with the metric
$H$ we obtain, 
$$
\aligned
&H^{\a\b} \pr_\a\pr_\b \pr H = \pr F + (\pr H^{\a\b}) 
\pr_\a\pr_\b \pr H = F^1 ,\\ 
&H^{\a\b} \pr_\a\pr_\b \pr G^\mu = \pr F_\mu +  (\pr H^{\a\b}) 
\pr_\a\pr_\b \pr G^\mu = F^1_\mu
\endaligned
$$
Using \eqref{2cuh3}-\eqref{2cuh4} and the inequality
$\|\pr H\|_{L^1_t L^\infty_x}\les \la^{-4\eps_0}$ of \eqref{ash1},
we infer that 
$$
\|F\|_{L^1_t L^2_x} \les \la^{-\f12},
\qquad \|F^1_\mu\|_{L^1_t L^2_x} \les \mu^{1-4\eps_0} \la^{-\f12-4\eps_0}
$$
Thus using the generalized energy identity for $\pr H$ and $\pr G^\mu$
we will have
$$
\il_{C_u}  |D_*\pr H|^2 \les
\|\pr^2 H\|^2_{L^\infty_t L^2_x} + \| \piT \|_{L^1_t L^\infty_x} 
\|\pr^2 H\|^2_{L^\infty_t L^2_x} +   \| F^1\|_{L^1_t L^2_x}
\|\pr^2 H\|_{L^\infty_t L^2_x}\les \la^{-1}
$$
Also,
$$
\aligned
\il_{C_u}  |D_*\pr G^\mu|^2 \les
\|\pr^2 G^\mu\|^2_{L^\infty_t L^2_x} + \| \piT \|_{L^1_t L^\infty_x} 
\|\pr^2 G^\mu\|^2_{L^\infty_t L^2_x} +  & \| F^1_\mu \|_{L^1_t L^2_x}
\|\pr^2 G^\mu\|_{L^\infty_t L^2_x}\\ &\les  \mu^{1-8\eps_0} \la^{-1-8\eps_0}
\endaligned
$$
\end{proof}

\section{A  remarkable property of $\rr_{44}$}

While the spacetime metric $\gg$  verifies the Einstein equations
$\rr_{\mu\nu}(\gg)=0$ this is certainly not true for the effective 
metric $H=H_{(\la)}$. This  could create serious  problems in the 
proof of the asymptotics theorem as the Ricci curvature appears
as a source term in the null structure equations. We have already established an
improved estimate for $\ric(H)$ in $L_t^1L_x^\infty$, see \eqref{ash6}.
This was done by comparing $\rr_{\mu\nu}(H)$ with $\rr_{\mu\nu}(G)=0$
where $G=\gg(\la^{-1}t,\la^{-1}x)$ is the rescaled Einstein metric. We need however
a stronger estimate involving the derivatives of $\rr_{44}(H)$
 along the null cones
$C_u$. To establish such an estimate we encounter an additional difficulty;
 the null cones $C_u$ have been constructed  relative to the approximate metric
  $H$.  This leads to significant differences between the $C_u$ energy estimates for
the second derivatives of $H$, see \eqref{cuH} and the corresponding ones\footnote{
The estimates for the second derivatives of the higher frequencies of  $G$
 do in fact diverge badly.} for $G$, see
\eqref{10cuh} in proposition \ref{CuH}. 
Using however the specific structure  of the component $\rr_{44}$
relative to the wave coordinates we can overcome this difficulty
and prove the following:
\begin{theorem} On  any null hypersurface $C_u$,
\be{ricci113}
\int_u^t\|\nab \rr_{44}(H)\|_{L^2(\stau)} d\tau\les \la^{-1}
\end{equation}
\label{Ricci4}
\end{theorem}
\begin{proof}  The proof of the theorem requires a rather
long and tedious argument which we present in our
paper \cite{Einst3}. 
\end{proof}

\section{Asymptotics theorem}
\label{se:AsThrm}
We start by recalling already established estimates for
the metric related quantities which play crucial role in what follows.
\begin{align}
&\|\pr H \|_{L^2_t L_x^\infty}\les
\la^{-\frac{1}2-4\ep_0},
\label{aso1}\\
&\|\pr H\|_{L^q(\stu)}\les \la^{\frac 2q-1-8(\frac 2q - \f12)\eps_0}\quad
\mbox{for}\quad  2\le q \le 4, \label{aso2}\\
&\|\ric(H)\|_{L^1_{t} L_x^\infty}\les 
\la^{-(1-4\ep_0)},
\label{aso3}\\
&\|\ric(H)\|_{L^p(\stu)}\les \la^{\frac 2p- 2- 8(\frac 2p - 1) \eps_0}\quad
\mbox{for}\quad  1\le p \le 2, \label{aso4}\\
&\|D_*\pr H\|_{L^2(C_u)}\les \la^{-\frac 12},\label{aso5}\\
& \il_0^s \|\nab \rr_{44}\|_{L^2(\stu)}\les \la^{-1-2\eps_0},\label{aso6}\\
&\|\rr_*\|_{L^2(C_u)}\les \la^{-\f12}\label{aso7}
\end{align} 
where
$\|\rr_*\|_{L^2(C_u)}:=\sum_{A,B,C,D}\|\rr_{ABCD}\|_{L^2(C_u)} +
\|\rr_{ABC4}\|_{L^2(C_u)}+
\|\rr_{B43A}\|_{L^2(C_u)}.
$
Note that some of the above estimates hold only throughout the region
$\Omega_*$.
\begin{theorem}
Throughout the region $\Omega_*$ the quantities $\trch -\frac 2r$, $\chih$,
and $\eta$
satisfy the following estimates:
\begin{equation}
\begin{align}
&\|\trch -\frac 2r\|_{L^2_t L^\infty_x}+ \|\chih\|_{L^2_t L^\infty_x} +
 \|\eta\|_{L^2_t L^\infty_x}\les \la^{-\frac 12-3\eps_0},\label{trih1}\\
&\|\trch -\frac 2r\|_{L^q(\stu)} + \|\chih\|_{L^q(\stu)} +
\|\eta \|_{L^q(\stu)}\les \la^{-3\eps_0}.
\label{trih2}
\end{align}
\end{equation}
In the estimate \eqref{trih1} function $\frac 2r$ can be replaced 
with $\frac 2{n(t-u)}$. We can also state the corresponding $L^1_t$
estimate following by H\"older inequality: 
\begin{equation}
\|\trch -\frac 2{n(t-u)}\|_{L^2_t L^\infty_x}\les \la^{-\frac 12-3\eps_0},
\qquad \|\trch -\frac 2{n(t-u)}\|_{L^1_t L^\infty_x}\les \la^{-3\eps_0}
\label{2trih4}
\end{equation}
In addition, in the exterior region $r\ge t/2$, 
\begin{equation}
\begin{split}
&\|\trch -\frac 2s\|_{L^\infty(\stu)}\les t^{-1} \la^{-4\eps_0},\qquad
\|\chih\|_{L^\infty(\stu)}\les t^{-1} \la^{-\eps_0} + \|\pr H(t)\|_{L^\infty_x},\\
&\|\eta\|_{L^\infty(\stu)}\les \la^{-1} + \la^{-\eps_0} t^{-1} 
+  \la^{\eps}\|\pr (H)(t)\|_{L^\infty_x}.
\end{split}
\label{trih3}
\end{equation}
where the last estimate holds for an arbitrary positive $\eps$, $\eps<\eps_0$. 
We also have the following estimates for the derivatives of $\trch$:
\begin{align}
& \|\sup_{r\ge \frac t2}
\|\Lb(\trch -\frac 2{r})\|_{L^2(\stu)}\|_{L^1_t} + \|\sup_{r\ge \frac t2}
\|\Lb(\trch - \frac 2{n(t-u)})\|_{L^2(\stu)}\|_{L^1_t}\les \la^{-3\eps_0 },
\label{2trih6}\\
&\|\sup_{r\ge \frac t2}\|\nabb \trch \|_{L^2(\stu)}\|_{L^1_t} +
\|\sup_{r\ge \frac t2}\|\nabb \big (\trch -\frac 2{n(t-u)}\big)
\|_{L^2(\stu)}\|_{L^1_t}\les \la^{-3\eps_0 }
\label{2trih8}
\end{align}

In addition we also have  weak  estimates of the form,
\be{lastz}
\sup_{u\le \frac t2}\|(\nabb, \Lb)\big( \trch -\frac
2{n(t-u)}\big)\|_{L^\infty(\stu)}\les \la^C
\end{equation}
for some large value of $C$.
\label{Trih}
\end{theorem}
\begin{corollary}
The estimates of theorem \ref{Trih} can be extended to the whole
region ${\mathcal I}_0^+\cap ([0,t_*]\times {\Bbb R}^3)$, where 
 ${\mathcal I}_0^+$ is the future domain of the origin on $\Si_0$.
\label{Ext}
\end{corollary}
\begin{remark}
The proof of the corollary \ref{Ext} requires an extension argument.
The estimates of the Asymptotics Theorem, which are uniform with respect
to the bootstrap region $\Omega_*$,  provide very good control 
of the foliations $C_u$ and $\stu$. By the standard continuity argument
this allows us to show that the estimates, in fact, hold in the maximal
domain allowed by the background estimates \eqref{aso1}-\eqref{aso7} on 
 the metric $H$, ${\mathcal I}_0^+\cap ([0,t_*]\times {\Bbb R}^3)$.
\end{remark}
\begin{proof} We shall first prove the estimates \eqref{trih1} and \eqref{trih2}
in the bootstrap region $\Om_{*}$. Once this is done we can easily
infer, by a standard continuity argument, that the estimates hold in fact in 
the entire region....

 To simplify our calculations we start
with  the following definition.
\begin{definition} We set,
\be{Thetadef}
\Theta=|\trch-\frac 2r| +|\trch-\frac 2s|+ |\chih| + |\eta| + |\pr H|
\end{equation}
\end{definition}
In view of our bootstrap assumptions ${\bf B_1)}$, 
 ${\bf B_2)}$ (  see section 6.1),  Remark \ref{RRs}  as well as the estimates 
\eqref{aso1}-\eqref{aso2} for $\pr H$ we  can freely make use of the
following:
\be{Thetaestim}
\|\Theta\|_{L_t^2L_x^\infty}\les \la^{-\f12 -2\ep_0} ,\qquad
 \|\Theta\|_{L^q(\stu)}\les \la^{-2\ep_0}
\end{equation} 
inside the bootstrap region $\Om_{*}$.

\subsection{ Estimates for $\trch, \chih$}\,\,\,
\label{se:trest}
\medn

We start with estimates \eqref{trih1}-\eqref{trih3} for $\tr\chi$. 
Observe that in view of the Corollary \ref{Comp}
it suffices to prove 
the desired estimates for $\trch -\frac 2s$.

Writing $y=(\tr\chi-\frac{2}{s})$
we have,
\be{8.98}
L(y) + \trch y=-R_{44} -\frac{2}{s}\bk_{NN} + \Theta^2
\end{equation}
Applying  the transport lemma \ref{Tran} 
we infer that at any point $P\in \Omega_*$,
$$
|s^2y(P)|\les \int_{\ga} s^2\bigg(|\rr_{44}|+\frac 1s |\pr H| +\Theta^2
\bigg) 
$$
where $\ga$ is the outgoing  null geodesic initiating on the time axis $\Ga_t$
 passing 
 through $P$ and $s$ is the  corresponding value of the affine parameter
$s$. Therefore,
$$| y(P)| \les\|\rr_{44}\|_{L_t^1L_x^\infty} + \frac{1}{s}
\int_{\ga}|\pr H|  +\|\Theta\|^2_{L_t^2 L_x^\infty}$$
and,  in view of  \eqref{Thetaestim}  and
\eqref{aso3},
\be{along1}
\|y(P)\|_{L^\infty}\les \la^{-1-4 \ep_0}+ \la^{-1-4\eps_0} + \frac{1}{s}
\int_{\ga}|\pr H|  
\end{equation}
In the exterior region $ s\ge \frac t2$, using the condition 
\eqref{aso1},
we infer that,
\be{2tr3}
\|\trch -\frac 2s\|_{L^\infty(\stu)}\les t^{-1} \la^{-4 \ep_0}.
\end{equation}
which proves \eqref{trih3}.
On the other hand, see  also the proof of  lemma \ref{Maxim},  \eqref{along1} 
leads a global estimate,
\be{2tr1}
\|\trch -\frac 2s\|_{L^\infty_x}\les \la^{-1-4 \ep_0} + \MH(t)
\end{equation}
where $\MH$ is  the maximal function of $\|\pr H(t)\|_{L_x^\infty}$.
The estimates \eqref{2tr1} and \eqref{aso1} together with
the corresponding  maximal function estimates readily imply that 
$$
\|\trch -\frac 2s\|_{L^2_t L^\infty_x} \les \la^{-\frac 12 -4\eps_0} + 
\| \MH(t)\|_{L^2_t} \les  \la^{-\frac 12 -4\eps_0} + \|\pr H\|_{L^2_t L^\infty_x}
\les \la^{-\frac 12 -4\eps_0}.
$$
On the other hand, using the comparison results between
$r$ and $s$, see section 6.3.,  $s\les \la^{1-8\ep_0}\les \la$, and
the H\"older inequalities
$$
\|\trch -\frac 2s\|_{L^q(\stu)} \les r^{\frac 2q} \|y\|_{L^\infty(\stu)} \les
\la^{\frac 2q} \la^{-1-4\eps_0} + s^{\frac 2q} s^{-\frac 12} \|\pr H\|_{L^2_t L^\infty_x}
\les \la^{\frac 2q} \la^{-1-4\eps_0}\les \la^{-4\eps_0}
$$
provided that $q> 2$ is chosen sufficiently close to $2$. 
Using the comparison results between $\frac{2}{r}$ and  $\frac{2}{s}$
of Corollary \ref{Comp}
we infer that,
\bea
\|\trch -\frac 2r\|_{L^2_t L^\infty_x}& \les &\la^{-\frac 12 -4\eps_0}\label{8.100}\\
\|\trch -\frac 2r\|_{L^q(\stu)} &\les&  \la^{-4\eps_0}\label{8.101}
\eea
as desired in \eqref{trih1} and \eqref{trih2}.
 Finally, \eqref{2trih4}
 follows from \eqref{stnu5} of Corollary
\ref{CoSn}.

We shall now estimate $\chih$ from the 
Codazzi equations \eqref{Codaz},
\be{8.99}
(\divv \chih)_A + \chih_{AB}k_{BN}=\half (\nabb_A \trch + k_{AN} \trch) -
\rr_{B{4}AB}.
\end{equation}
Taking advantage of corollary \ref{corB4AB}, with a different error
 term $E$,  we rewrite it  in
the form,
\be{CZZ}
(\divv \chih)_A =\half \nabb_A(\trch-\frac{2}{r}) +\nabb_A\pi+\nabb^B\pi_{AB} + E
\end{equation}
with $\pi$ and $E$ obeying pointwise estimates 
$$
|\pi |\les |\pr H|,\qquad E\les  \Theta\cdot \pr H + \frac 1r |\pr H|
$$
We shall now take advantage of the elliptic 
estimate of
proposition \ref{CZestimates}; we 
write,
\bea
\| \chih\|_{L^\infty(\stu)}&\les&\la^\ep
\|\trch-\frac{2}{r}\|_{L^\infty(\stu)}\nn\\
&+&\la^\ep\|\pi\|_{L^\infty(\stu)}
+r^{1-\frac 2 q}\|E\|_{L^q(\stu)}\label{eq:200}
\eea
with $q>2$. 
\begin{remark}
In the application of the elliptic estimate \eqref{caldzygm2} in the 
derivation of \eqref{eq:200} we need some rough estimates 
for $\nabb\trch$ of the type 
$$
\|r\nabb\trch\|_{L^\infty(\stu)}\les \la^{C}
$$
for some large constant $C>0$.
These  weak estimates, consistent with \eqref{lastz},
 are a lot easier to derive and 
can be obtained directly from the transport equations 
\eqref{D4trchi}, \eqref{D4chih} for $\trch$ and $\chih$. 
We refer the
reader to our paper \cite{Kl-Ro} for more details.
\label{1weakremark}
\end{remark}

Therefore, choosing $q=2+\ep\,$ for  sufficiently small $\ep>0$, and
using the bootstrap assumptions ${\bf B2)}$  as well as the assumptions 
\eqref{aso2}
we infer that,
\begin{equation}
\begin{split}
\| \chih\|_{L^\infty(\stu)}&\les
\la^{\ep}\|\trch-\frac{2}{r}\|_{L^\infty(\stu)}
+\la^{\ep}\|\pr H \|_{L^\infty_x}
\\
&+r^{1-\frac 2 q}\bigg( \|\Theta\|_{L^q(\stu)}
\|\pr H\|_{L^\infty_x} + r^{-1+\frac 2 q}\|\pr H\|_{L^\infty_x} \bigg)\\
&\les \la^{\ep}\bigg (\|\trch-\frac{2}{r}\|_{L^\infty(\stu)}
+ \|\pr H \|_{L^\infty_x}\bigg ) 
\end{split}
\label{2ch3}
\end{equation}

Now we observe that the desired pointwise estimate \eqref{trih3} in the exterior
region
$r\ge t/2$ follows from \eqref{2tr3} and the estimate 
$|\frac 2r-\frac 2s|\les \la^{-\eps_0}s^{-1}\les \la^{-\eps_0}t^{-1}$,
\be{2ch2}
\| \chih\|_{L^\infty(\stu)}\les t^{-1} \la^{-\eps_0} + \|\pr H\|_{L^\infty_x}
\end{equation}
We can also add a global estimate following from Corollary \ref{Comp}\footnote
{Namely, the inequality $\|\trch -\frac 2r\|_{L^\infty_x} \les  {\cal M}^3
\big(\|\trch -\frac 2s\|_{L^\infty_x}\big )$} and
\eqref{2tr1}.
\be{2ch4}
\| \chih\|_{L^\infty(\stu)}\les \la^{-1-4\eps_0} + \pr H(t)+
{\cal M}^4(\pr H) (t).
\end{equation}
Now squaring and integrating \eqref{2ch3} in time we infer from \eqref{aso1}
and the just proved estimate \eqref{8.100} for $\trch -\frac 2r$ that
\be{2ch1}
\| \chih\|_{L^2_t L^\infty_x}\les \la^\eps \bigg
(\|\trch-\frac 2r\|_{L^2_t L^\infty_x} +  \|\pr H \|_{L^2_t L^\infty_x}\bigg ) 
\les \la^{-\frac 12-3\eps_0},
\end{equation}
 which is the estimate claimed in \eqref{trih1} of theorem\ref{Trih}.

On the other hand, applying the elliptic estimate \eqref{caldzygm1} of proposition
\ref{CZestimates} to the equation \eqref{CZZ} yields the following:
$$
\|\chih\|_{L^q(\stu)}\les \|\trch -\frac 2r\|_{L^q(\stu)} + 
\|\pr H\|_{L^q(\stu)} + \|E\|_{L^p(\stu)}
$$
for some $q\ge 2$, $\frac 1p=\frac 12+\frac 1q$.
Choosing $q=2+\eps\,$ as in bootstrap assumption ${\bf B2)}$ we infer 
with the help of the  estimate \eqref{8.101} for $\trch - \frac 2r$ 
and \eqref{aso2},
that
$$
\aligned
\|\chih\|_{L^q(\stu)}&\les \la^{-4\eps_0} 
+  \|\Theta \,\pr H\|_{L^{p}(\stu)} +  \frac 1r \|\pr H\|_{L^{p}(\stu)}\\
 &\les
\la^{-4\eps_0} + \|\pr H\|_{L^{2}(\stu)}\|\Theta\|_{L^q(\stu)} 
 + \|\pr H\|_{L^{q}(\stu)}\\ &\les 
\la^{-4\eps_0}.
\endaligned
$$

\subsection{ Estimates for $\eta$}\,\,\,\,
\medn

 We start  with the Hodge system \eqref{diveta}--\eqref{curleta}:
\beaa
\divv\,\eta &=& \half\bigg(\mu +2\bk_{NN}\trch   -2|\eta|^2 -|\chih|^2
-2k_{AB}\chi_{AB}\bigg)
 - \half \de^{AB}\rr_{A{43}B},\\
\curll\,\eta &=& \half\in^{AB} k_{AC} \chih_{CB} -
 \half \in^{AB}\rr_{A{43}B}
\eeaa
with $\mu$ defined in \eqref{eqmu},
$\mu=\Lb(\trch)-\f12 (\trch)^2 - 
\big (k_{NN}+n^{-1}\nab_N n\big )\trch$
and satisfying the transport equation \eqref{D4tmu},
\begin{equation}
\begin{split}
L(\mu) + \trch \mu &=2(\etab_A-\eta_A)\nabb_A(\trch)
-2\chih_{AB}\Bigl (2\nabb_A \eta_B +  2\eta_A\eta_B   \\&+ \bk_{NN}\chih_{AB} 
+\trch\chih_{AB} +\chih_{AC}\chih_{CB} +2k_{AC}\chi_{CB}+ \rr_{B{43}A}\Bigr )
\\&
-\Lb(\rr_{44})+  (2k_{NN} - 4 n^{-1}\nab_N n)
)\big (\half (\trch)^2 - |\chih|^2 \\ &-
\bk_{NN} \trch -\rr_{44}\big ) +4\bk_{NN}^2\trch 
+(\trch+4\bk_{NN})(|\chih|^2+\rr_{44}) \\ &- \trch 
\bigg (2 (k_{AN} - \eta_A) n^{-1} \nab_A n - 2 |n^{-1} N(n)|^2 +
\half \rr_{4343} + 2 k_{Nm} k^{m}_N\bigg)
\end{split}
\label{2eta1}
\end{equation}
Observe that in view of Corollary \ref{corA4B3}
we can rewrite our div-curl system for $\eta$ as follows:
\bea
\divv\,\eta &=&\divv\pi^{(1)} +  \half\bigg(\mu +2\bk_{NN}\trch 
  -2|\eta|^2 -|\chih|^2
-2k_{AB}\chi_{AB}\bigg)
 - \frac 12\rb + E^{(1)}, \nn\\
\curll\,\eta &=&\curll\pi^{(2)}+ \half\in^{AB} k_{AC} \chih_{CB} + E^{(2)}.
\label{8.103}
\eea
where $\rb=(\rr+\rr_{34})$ and  
$$
\aligned
&|\pi^{(1,2)}|\les |\pr H|\\
&|E^{(1,2)}|\les (|\pr H|^2+|\chi||\pr H|).
\endaligned
$$
\begin{remark} We would like  to treat the system  formed by the transport
 equation \eqref{2eta1} coupled with the elliptic system \eqref{8.103}
in the same manner as we have dealt with the system for $\trch$ and $\chih$.
Indeed the Hodge system \eqref{8.103} is similar to the Hodge system \eqref{8.99}.
The transport equation for $\mu$ differs however significantly from the
transport equation \eqref{8.98} for $\trch$. Indeed the  only curvature term
on the right hand side of \eqref{8.98} is  $\rr_{44}$ while the right hand side
of \eqref{2eta1} exhibits the far more dangerous term.
$\Lb(\rr_{44})$. In what follows we shall get around this difficulty by introducing a
new covector $\muv$ through a Hodge system on the surfaces $S_{t,u}$.
Using once more the special structure of the Einstein equations we shall
derive a new  transport equation for $\muv$ whose right hand side 
exhibits only terms depending on $\ric(H)$ and favorable components 
of the curvature tensor.
\end{remark}
We define an auxiliary $S$-tangent co-vector $\muv_A $ as a solution of
the Hodge system 
\bea
\divv\muv &=& \mu-\rb,\label{musla1}\\
\curll\muv &=& 0 \label{musla2}
\end{eqnarray}
with $\rb=\rr_{43} +\rr$.
We now prove the following 
\begin{proposition}
\medn
\begin{enumerate}
\item
The covector $\muv$ verifies the following,
\beaa
\divv \big(\ddd_4\muv + \half \trch\muv -\chih\cdot \muv \big ) &= &     
\pr H\cdot \ddd_{4}\muv + \nabb_A\bigg(2\rr_{A4}
+ \frac{2}r\pi_A +  \Theta\cdot \Theta\bigg) \\&-& 
\frac 2r( 3\rr_{34}  +2\rr)  
+\Theta \ric +\Theta \rr_* + \Theta\cdot D_*\pr H \nn \\
& +& \Theta\cdot\Theta\cdot \Theta 
 + \frac 1r \Theta\cdot\Theta + \frac 1{r^2}\pr H ,\nn\\
\curll \big(\ddd_4\muv + \half \trch\muv -\chih\cdot \muv \big ) &=&
\pr H\cdot \ddd_{4}\muv + \nabb \big (\Theta\cdot \Theta  \big) +
\rr_*\cdot\Theta \\ &+& \frac{1}{r}\Theta\cdot\Theta +\Theta \cdot \Theta \cdot 
\Theta \nn
\eeaa
\item The covector $\muv$ verifies the following estimates
\be{10mu}
\|\muv\|_{L^\infty(\stu)}\les \la^{-1} + \cal M(\pr H) 
\end{equation}
\be{10mu1}
\|\muv\|_{L^q(\stu)} \les\la^{-3\eps_0}
\end{equation}
\end{enumerate}
\label{funnymu}
\end{proposition}

\textbf{Proof of part 2 of proposition \ref{funnymu}}\,\,
\begin{remark} For convenience we extend our bootstrap
assumptions ${\bf B1)}$ and ${\bf B2)}$ to include $\muv$.
Thus, throughout the proof below,  we redefine $\Theta$, see \eqref{Thetadef}, 
as follows:
\be{Thetadefnew}
\Theta=|\trch-\frac 2r| +|\trch-\frac 2s|+ |\chih| + |\eta| + |\pr H| +|\muv|
\end{equation}
This is justified  since our stated  estimates  are stronger
than ${\bf B1)}$ and ${\bf B2)}$ for $\muv$.
\end{remark}
 Assuming
the first part of the proposition \ref{funnymu} we now derive 
the estimates of part 2. 
We start by applying  the elliptic estimates of 
proposition \ref{CZestimates} to the Hodge system
of proposition \ref{funnymu}. Thus
for some $q>2$, denoting  by $M$
the quantity 
$$M=\big(\ddd_4\muv + \half \trch\muv -
\chih\cdot \muv\big ),$$
we have,
$$
\aligned
\|M\|_{L^\infty(\stu)}&\les \|\pr H\|_{L^q(\stu)} \|M\|_{L^\infty(\stu)} +
\la^\eps \bigg(\|\ric(H)\|_{L^\infty(\stu)}
+\|\Theta\|^2_{L^\infty(\stu)} +
\frac 1r \|\pr H\|_{L^\infty(\stu)}\bigg)\\
& +
r^{1-\frac 2q} \bigg(\|\Theta \rr_*\|_{L^q(\stu)} +
\|\Theta \nabb (\pr H)\|_{L^q(\stu)}
 +  
\|\Theta \ric(H)\|_{L^q(\stu)}\\
& + \frac 1r\|\ric(H)\|_{L^q(\stu)}
 +\|\Theta^3\|_{L^q(\stu)} + \frac 1r\|\Theta^2\|_{L^q(\stu)} + 
\frac 1{r^2}\|\pr H\|_{L^q(\stu)}\bigg).
\endaligned
$$
\begin{remark}
As in the case of the estimates for $\chih$, the use of the elliptic
estimates \eqref{caldzygm2} of proposition \ref{CZestimates} for the 
Hodge system satisfies by the quantity $M$ requires rough 
estimates of the type 
$$
\|r \nabb \rr_{A4}\|_{L^\infty(\stu)} + \|\nabb \pi\|_{L^\infty(\stu)}
+ \|r\Theta \cdot \nabb\Theta\|_{L^q(\stu)}\les  \la^{C}
$$
for some $q>2$.  The estimate for the derivatives of the 
Ricci curvature and the metric $H$ are contained in our background
estimates \eqref{ash1}-\eqref{ash6}. In addition to $\trch$ and $\chih$,
for which we have already outlined the procedure of obtaining 
such weak estimates, the quantity $\Theta$ contains $\eta$ and $\muv$.
Once again, we can use the transport equation \eqref{D4eta} for $\eta$
and the Hodge system \eqref{musla1}-\eqref{musla2} combined with the transport equation 
\eqref{2eta1} for $\mu$  to handle these terms.
\end{remark}
 
 Taking  $q$  
sufficiently close to $q=2$,  using the bootstrap assumption, 
 $\|\Theta\|_{L^q(\stu)}\les \la^{-2\eps_0}\les 1$, and the estimate 
 $\|\pr H\|_{L^q(\stu)}\les \la^{-2\eps_0}< 1/2$
 we can then conclude that
$$
\aligned
\|M\|_{L^\infty(\stu)}&\les  \la^{\eps}
\bigg(\|\ric(H)\|_{L^\infty(\stu)} +
\|\Theta\|^2_{L^\infty(\stu)}  
+\frac 1r \|\pr H\|_{L^\infty(\stu)}\\&
+ \|\Theta\|_{L^\infty(\stu)}\|D_*\pr H\|_{L^q(\stu)} + 
\|\Theta\|_{L^\infty(\stu)}\|\rr_*\|_{L^q(\stu)}\bigg )
\endaligned
$$

Applying  the transport lemma \ref{Tran} to the transport equation
\be{8.110}
\ddd_4\muv + \half \trch\muv =M+
\chih\cdot \muv,
\end{equation}
we infer that at any point $P\in \Omega_*$,
$$
|s\muv(P)|\les \int_{\ga} s\bigg(|M| +\Theta^2\bigg) 
$$
where $\ga$ is the outgoing  null geodesic initiating on the time axis $\Ga_t$
 passing 
 through $P$ and $s$ is the  corresponding value of the affine parameter
$s$.
Hence, 
\beaa
|\muv(P)|&\les&\la^{\eps}\int_\ga
\bigg( 
\|\Theta\|^2_{L^\infty(\stu)}  
+ \|\Theta\|_{L^\infty(\stu)} \|D_*\pr H\|_{L^q(\stu)}\\ &+& 
\|\Theta\|_{L^\infty(\stu)}\|\rr_*\|_{L^q(\stu)}\bigg ) +
\frac 1s\int_\ga \|\pr H\|_{L^\infty(\stu)}
\eeaa
Observe that by ${\bf B1)}$ and \eqref{aso5}, we have
$$
\aligned
\il_{\ga} \|\Theta\|_{L^\infty(\stu)} \|D_*\pr H\|_{L^q(\stu)}&\les
 \|\Theta\|_{L^2_t L^\infty_x}\bigg(\il_{\ga} 
\|D_*\pr H\|^2_{L^q(\stu)}\bigg)^{\half}\\
&\les \la^{-\f12 - 2\eps_0}\|\pr^2 H\|_{L^2_t L_x^\infty}^{1-\frac 2q} 
\|D_*\pr H\|_{L^2(C_u)}^{\frac 2q} \\ & \les \la^{-\f12- 2\eps_0}
\la^{-(1+ 4\eps_0)(\frac 12 -\frac 1q)}\la^{-\frac 1q}\\
&\les \la^{-1-2\eps_0}
\endaligned
$$ 
A similar estimate, by \eqref{aso7}, also holds for the term 
involving $\rr_*$. 
Consequently,
$$
\|\muv\|_{L^\infty(\stu)}\les \la^{-1} + \cal M(\pr H) 
$$
as desired.

Observe also that in the exterior region $r\ge \frac{t}{2}$,
\be{muvext}
\|\muv\|_{L^\infty(\stu)}\les \la^{-1} + r^{-1}\la^{-4\ep_0}. 
\end{equation}

Going back to proposition \ref{funnymu} and 
applying now the estimate \eqref{caldzygm1}
of proposition \ref{CZestimates} we deduce,
 for $\frac 1p=\frac 12+\frac 1q$,
$$
\aligned
\|M\|_{L^q(\stu)}&\les \|\pr H\|_{L^2(\stu)}\|M\|_{L^q(\stu)}+ 
\|\ric(H)\|_{L^{q}(\stu)} +\frac 1r  \|\pr H\|_{L^q(\stu)}+
\|\Theta\|^2_{L^{2q}(\stu)}\\ &  
+ \|\Theta\|_{L^q(\stu)} \bigg(\|D_*\pr H\|_{L^2(\stu)}+
\|\rr_*\|_{L^2(\stu)}\bigg) \\
&+ r^{\frac 2p-1} \|\Theta\|^2_{L^\infty(\stu)} + 
r^{\frac 2p -2}\|\pr H\|_{L^\infty(\stu)}
\endaligned
$$
According to the estimates \eqref{aso2},
$
\|\pr H\|_{L^2(\stu)}\les  \la^{-4\eps_0}<1
$. Thus we can absorb the term with $M$ into the left hand-side. 

On the other hand,  using the transport lemma \ref{Tran}
applied to the transport  equation \eqref{8.110} we infer,

$$\|\muv\|_{L^q(\stu)}\les \frac 1{r(t)^{(1-\frac 2q)}}\il_u^t 
r(t')^{(1-\frac 2q)}\bigg(
\|M\|_{L^q(S_{t', u})}+ \|\Theta\|^2_{L^{2q}(S_{t', u})} \bigg) dt'
$$

Applying
the  bootstrap assumptions
${\bf B1)}$, ${\bf B2)}$, and \eqref{aso1}-\eqref{aso7} we infer
that,
\begin{equation}
\begin{split}
\|\muv\|_{L^q(\stu)} &\les \la^\eps\bigg(
\la^{\f12}\|\ric( H)\|^{\half}_{L^1_t L^\infty_x}\|\ric( H)\|^{\half}_{L^{\frac q2}(\stu)} + 
\|\Theta\|_{L^1_t L^{\infty}_x}\|\Theta\|_{L^{q}(\stu)}\\
&+\frac 1{r^{(1-\frac 2q)}}\il_u^t r(t')^{(1-\frac 2q)}r(t')^{\frac 2q -1}
\|\pr H\|_{ L^\infty_x(S_{t',u})}dt'\nn\bigg)\\
 &+ r^{\half}\|\Theta\|_{L^{q}(\stu)}\|D_*\pr H\|_{L^2(C_u)}
+ r^{1-\frac {2p}q}\|\Theta\|^2_{L^{q}(\stu)}\|\Theta\|_{L^1_t L^{\infty}_x} 
\\ &+r^{\frac 2p-1} \|\Theta\|^2_{L^2_t L^\infty_x} + 
\frac 1{r^{(1-\frac 2q)}}\il_u^t r(t')^{(1-\frac 2q)}r(t')^{\frac 2p -2}
\|\pr H\|_{ L^\infty_x(S_{t',u})}dt'\nn\\ &\les 
\la^{-3\eps_0} + r^{\frac 2p - \frac 32}\|\pr H\|_{L^2_t L^\infty_x}\les
\la^{-3\eps_0}
\end{split}
\end{equation}
as desired.
On the right hand-side of the last series of inequalities, for the sake of brevity, 
we have abused the notation using $\|f\|_{L^q(\stu)}$ to denote 
$\sup_{t,u} \|f\|_{L^q(\stu)}$.

\bigskip
Using the estimates \eqref{10mu} and \eqref{10mu1} for $\muv$ 
we are now ready to return to the proof of the estimates for $\eta$

Now with the help of the established  estimates for  $\muv$ we
 shall derive the desired
estimates for 
$\|\eta\|_{L^2_t L^\infty_x}$ and $\|\eta\|_{L^q(\stu)}$.
First observe that using using the definition \eqref{musla1}-\eqref{musla2}
of $\muv$ the div-curl system \eqref{8.103} for $\eta$ takes the form
\beaa
\divv\,(\eta - \frac 12 \muv)&=&\divv\pi^{(1)} +  
\frac 1r \pr H   + \Theta\cdot \Theta ,\\
\curll\,(\eta -\frac 12 \muv) &=&\curll\pi^{(2)}+ \frac 1r \pr H   + \Theta\cdot \Theta
\eeaa
We are now ready to apply proposition \ref{CZestimates}
to our Hodge system for $\eta - \frac 12 \muv$. Thus, for some $q>2$,
sufficiently close to $2$,
\beaa
\|\eta - \frac 12 \muv\|_{L^\infty(\stu)}&\les &
\la^\eps \|\pr H\|_{L^\infty(\stu)} + r^{-\frac 2q}\|\pr H\|_{L^q(\stu)} 
+r^{1-\frac 2 q}
\|\Theta^2\|_{L^q(\stu)}\\
&\les &\la^{\eps} \|\pr H\|_{L^\infty(\stu)} + 
\la^{-\eps_0}\|\Theta\|_{L^\infty(\stu)},
\eeaa
where we have used the bootstrap estimate $\|\Theta\|_{L^q(\stu)}\les \la^{-2\eps_0}$.
Furthermore,  we infer with the 
help of \eqref{10mu} that 
\begin{equation}
\begin{split}
\|\eta \|_{L^\infty(\stu)}&\les \la^{-1} + {\cal M(\pr H)} 
+ \la^{\eps} \|\pr H\|_{L^\infty(\stu)} +\la^{-\eps_0} \|\Theta\|_{L^\infty(\stu)}
\end{split}
\label{10et}
\end{equation}
The desired $L_t^2L^\infty_x$ estimate follows immediately
from the bootstrap assumption ${\bf B1)}$ and the estimates
\eqref{aso1}-\eqref{aso7}.

Consider also the exterior region $r\ge \frac{t}{2}$
Observe  that, using the estimates \eqref{2tr3}, \eqref{2ch2} and \eqref{muvext}  for 
$\trch-\frac{2}{r}, \chih, \muv$
 already established in the exterior region 
we infer that,

\begin{equation}
\|\eta \|_{L^\infty(\stu)}\les \la^{-1} + \la^{-\ep_0}t^{-1} 
+ \la^{\eps} \|\pr H\|_{L^\infty_x}
\label{10et'}
\end{equation}

\medn
On the other hand, for $\frac 1p=\half+\frac 1q$,
\beaa
\|\eta -\frac 12 \muv\|_{L^q(\stu)}&\les &
\la^\eps\|\pr H\|_{L^q(\stu)}+ \frac{2}{r}\|\pr H\|_{L^p(\stu)}\\
+\|\Theta^2\|_{L^p(\stu)}
\eeaa
Since $\frac 1p = \frac 12 + \frac 1q$ and $q\ge 2$, we have
$2p\le q$ and the H\"older inequality gives 
$$
\|\Theta \|_{L^{2p}(\stu)}^2\les r^{2-\frac {4}q} \|\Theta\|_{L^{q}(\stu)}^2
\les \la^{-3 \eps_0}
$$
from the bootstrap assumption ${\bf B2)}$, provided that $q$ is sufficiently close to 2. 
Thus with the help of \eqref{aso2} and  the estimate \eqref{10mu1} we obtain,
\begin{equation}
\|\eta\|_{L^q(\stu)}\les \|\muv\|_{L^q(\stu)} + \la^{-3 \eps_0}\les \la^{-3\eps_0}
\label{10et1}
\end{equation}
as desired.

{\bf Proof of part 1 of  proposition \ref{funnymu}}\,\,\,
We now concentrate on the proof of proposition \ref{funnymu}.
We start by expressing
the transport equation \eqref{2eta1} for $\mu=\Lb(\trch)-\f12 (\trch)^2 - 
\big (k_{NN}+n^{-1}\nab_N n\big )\trch$
in a more tractable form.
The troublesome terms are  $\Lb \rr_{44}$ and
$\trch \rr_{4343}$.
We shall first eliminate $\Lb R_{44}$ in exchange 
for more favorable terms.
We do this with the help of the twice contracted  Bianchi
identity:
$$D^\nu(\rr_{\mu\nu}-\frac{1}{2} g_{\mu\nu} \rr)=0$$
with $\rr$ the scalar curvature $\rr=g^{\mu\nu}\rr_{\mu\nu}$.
Thus, relative to our canonical null frame,
$$D^3\rr_{43}+D^4\rr_{44}+D^A\rr_{4A}=\half L(\rr),$$
or,
$$D_3\rr_{44}=-D_4\rr_{43}+2D^A\rr_{4A}- L(\rr).$$
On the other hand,
\beaa
D_3\rr_{44}&=&\Lb \rr_{44}-4\eta_A\rr_{A4}-2\bk_{NN} \rr_{44}\\
D_4\rr_{43}&=&L \rr_{43}-2\etab_A\rr_{4A}\\
D^A\rr_{4A}&=&\nabb^A\rr_{4A}-\chi_{AC}\rr_{CA}+k_{AN}\rr_{4A}-\half \trch\rr_{43}-
\half\trchb\rr_{44}
\eeaa
Therefore,
\beaa
\Lb(\rr_{44})&=&-L(\rr_{43}+\rr)+2\nabb^A\rr_{4A}\\&-&(2\rr_{34}+\rr-\rr_{44})\cdot \trch
+\ric\cdot(\chih, k,\eta)
\eeaa
Using this formula we can rewrite the transport equation for 
$\mu$ in the form:
\beaa
L(\mu) + \trch \mu = L(\rb) + 2 \nabb_A \rr_{A4} +\trch( 2\rr_{34} +\rr)
 -  \trch \rr_{4343}&&\nn \\
+ 2 (\etab_A - \eta_A)\nabb_A \trch -4 \chih\cdot \nabb \eta  +\Theta \cdot \rr_*
+ \Theta\cdot\Theta\cdot \Theta  +\Theta \cdot\ric+ \frac 1r \Theta\cdot\Theta + \frac
1{r^2}\pr H,&&
\eeaa
where $\rb = \rr_{43} +\rr$.
Thus
\bea
L(\mu -\rb) + \trch (\mu -\rb)= 2 \nabb_A \rr_{A4} +\trch \rr_{34}  - 
 \trch \rr_{4343}+  \Theta \cdot\ric&&\label{eqmr} \\
+
2 (\etab_A - \eta_A)\nabb_A \trch -4 \chih\cdot \nabb \eta +\Theta \cdot \rr_*
+ \Theta\cdot\Theta\cdot \Theta  + \frac 1r \Theta\cdot\Theta + \frac 1{r^2}\pr H&&
\nn 
\end{eqnarray}
Observe that $\rr_{34}=H^{\a\b}\rr_{\a 3 \b 4}=\f12 \rr_{4343}-\de^{AB}\rr_{A34B}$.
Also, $\rr_{AB}=-\f12 \rr_{3A4B}-\f12 \rr_{4A3B}+\de^{CD}\rr_{CADB}$.
Therefore,
$$\rr_{3434}=2(\rr_{34}+\de^{AB}\rr_{AB})+\de^{CD}\de^{AB}\rr_{CADB}$$
or, since $\de^{AB}\rr_{AB}=\rr_{34}+\rr$, using 
corollary \ref{CorABCD} for $\rr_{ABCD}$,
$$\rr_{3434}=2(2\rr_{34}+\rr)-\divv \pi +E,$$
where  $$|\pi|\les|\pr H|\qquad  \mbox{and} \qquad |E|\les |\pr H|^2 +|\chi||\pr H|.$$
Using this we can rewrite \eqref{eqmr} in the form,
\bea
L(\mu -\rb) + \trch (\mu -\rb)&=& 2 \nabb_A \rr_{A4}+\trch \divv\pi 
-\trch( 3\rr_{34}  +2\rr)  
 \label{eqmr2} \\
&+&  2 (\etab_A - \eta_A)\nabb_A \trch -4 \chih\cdot \nabb \eta +
\Theta \cdot\ric+\Theta \cdot \rr_*
\nn\\
& +&  \Theta\cdot\Theta\cdot \Theta + \frac 1r \Theta\cdot\Theta + \frac 1{r^2}\pr H
\nn 
\end{eqnarray}

Recall that we defined an $S$-tangent co-vector $\muv_A $ as a solution of
the Hodge system
\bea
\divv\muv &=& \mu-\rb,\\
\curll\muv &=& 0
\end{eqnarray}
We shall now use the commutation formula of  lemma \ref{2Comm}.
\beaa
\divv(\ddd_4\muv)-L(\divv \muv)&=&\f12 \trch \divv\muv +\chih\cdot\nabb \muv-
n^{-1}\nabb n \cdot \ddd_4\muv\\
&+&\f12\trch \bk_{AN}\muv_A-\chih_{AB} \bk_{BN}\muv +\rr_{AB4B}\muv_A\\
\curll(\ddd_4\muv)-L(\curll \muv)&=&\f12 \trch \curll\muv +\in^{BA}\chih_{BC}\nabb_C\muv_A
-\in^{BA} n^{-1}\nabb_B n \ddd_4\muv_A \\&-& \in^{BC}\chi_{BA} \bk_{CN} \muv_A + 
\in^{BC} \rr_{AC4B}\muv_A
\eeaa
Using the transport equation \eqref{eqmr} and
commuting $L$ with $\divv$ and $\curll$ (see lemma \ref{2Comm})
we can derive the following Hodge system
for $\ddd_4 (\muv)$:
\begin{eqnarray}
\divv (\ddd_4\muv ) &=& -\f12 \trch \,\divv\muv + \chih\cdot \nabb\muv +
\pr H\cdot \ddd_{4}\muv + \frac 2r \divv\pi\nn \\
&+&2 (\etab_A - \eta_A)\nabb_A \trch -4 \chih\cdot \nabb \eta   + 
\rr_{AB4B}\muv_A  + 2\nabb_A \rr_{A4}\nn \\ &-&\frac 2r( 3\rr_{34}  +2\rr)  
+\Theta \ric +\Theta \rr_* + \Theta \cdot D_*\pr H
\nn \\
& +& \Theta\cdot\Theta\cdot \Theta 
 + \frac 1r \Theta\cdot\Theta + \frac 1{r^2}\pr H,\nn \\
\curll (\ddd_4\muv) &=& \half \trch\curll\muv \,+ \pr H\cdot \ddd_{4}\muv +
\in^{BA}\nabb_C\mu_A\chih_{BC}
\nn \\ &+& \in^{CB}\rr_{AB4C}\muv_A 
+\frac 1r \Theta\cdot \Theta + \Theta\cdot\Theta\cdot \Theta \nn 
\end{eqnarray}
\begin{remark}
We got rid of  the dangerous term $\Lb(\rr_{44})$. We still need 
to eliminate the terms of the form $\Theta \cdot \nabb\Theta$.
\end{remark}
Observe that according to the Codazzi equation
$$
\divv\chih_A - \f12\nabb_A\trch = \half k_{AN} \trch - \chih_{BN} k_{BN} - \rr_{B4AB}.
$$
Therefore,
\beaa
-\f12 \trch \,\divv\muv + \chih_{AB}\nabb_B
\muv_A &=&-\f12 \divv(\trch\muv)+\nabb^B(\chih_{AB}\muv_A)+\f12 \muv\cdot \nabb\trch-
(\nabb^B\chih_{AB})\muv_A\\
&=&-\nabb^A\bigg(\f12 \trch\muv_A-\chih_{AB}\muv_B\bigg)-\muv\bigg(\half k_{AN} \trch -
\chih_{BN} k_{BN} - \rr_{B4AB}\bigg)
\eeaa
Thus
\begin{eqnarray}
\divv (\ddd_4\muv + \half \trch\muv_A - \chih_{AB}\muv_B ) &=&  
\pr H\cdot \ddd_{4}\muv + \frac 2r \divv\pi
+2 (\etab_A - \eta_A)\nabb_A \trch -4 \chih\cdot \nabb \eta \nn\\ & +& 
 2\nabb_A \rr_{A4}- \frac 2r( 3\rr_{34}  +2\rr)  
+\Theta \ric +\Theta \rr_*\nn \\
& +& \Theta\cdot D_*\pr H + \Theta\cdot\Theta\cdot \Theta 
 + \frac 1r \Theta\cdot\Theta + \frac 1{r^2}\pr H,\label{optical1}
\end{eqnarray}

Also, since $\curll \muv =0$, we have
$$
\aligned
\half \trch\curll\muv &+\in^{BA}\nabb_C\muv_A\chih_{BC} = 
\in^{BA}\nabb_A\muv_C\chih_{BC} \\ &=  -\in^{AB} \nabb_A ( 
\chih_{BC} \muv_C)  + \in^{AB}  \nabb_A \chih_{BC}\muv_C \\ &=
-\in^{AB} \nabb_A (\half \trch \muv_B + \chih_{BC} \muv_C) +
\in^{AB} \nabb_A (\trch )  \muv_B   + \in^{AB}  (\nabb_A \chi_{BC}) \muv_C
\endaligned
$$ 
On the other hand,  see \cite{Kl-Ro} section 2,
$\nabb_A\chi_{BC}=\nabb_C\chi_{AB}-\rr_{B4CA}+k\cdot\chi$. Therefore,
\bea
\curll (\ddd_4\muv + \half \trch\muv_A - \chih_{AB}\muv_B )&=&
\pr H\cdot \ddd_{4}\muv - 2 \curll (\chih \cdot \muv)+ \in^{AB} \nabb_A (\trch )  \muv_B
\nn\\ &+& 
+\frac 1r \Theta\cdot \Theta + \Theta\cdot\Theta\cdot \Theta\label{optical2} 
\eea
Observe also, in \eqref{optical1}, using Codazzi
\beaa
  - 2\eta_A\nabb_A \trch -4 \chih^{AB}\cdot \nabb_B
\eta_A&=&- 2\eta_A\nabb_A \trch-4\nabb^A(\chih_{AB}\eta_B) +
4\nabb^A\chih_{AB}\eta_B\\
& =&-4\nabb^A(\chih_{AB}\eta_B)+4\eta_B\rr_{A4BA}+\eta\cdot \chi\cdot k
\eeaa 
Therefore,
\beaa
\divv \big(\ddd_4\muv + \half \trch\muv_A -\chih_{AB}\muv_B \big)& = &   
\pr H\cdot \ddd_{4}\muv +  
\nabb_A\bigg(2\rr_{A4}- 4\chih_{AB}\eta_B +\frac 2r \pi_A\bigg) \nn \\
 &+& 2\etab_A \nabb_A \trch - 
\frac 2r( 3\rr_{34}  +2\rr)  
+\Theta \ric +\Theta \rr_*  \nn \\
& +&  \Theta\cdot D_*\pr H + \Theta\cdot\Theta\cdot \Theta 
 + \frac 1r \Theta\cdot\Theta + \frac 1{r^2}\pr H\nn
\eeaa
In addition, since $\etab_A = -\bk_{AN}$, 
$$
\aligned
\etab_A \nabb_A\trch &= -\nabb_A \bigg(\bk_{AN} (\trch -\frac 2r)\bigg) + 
(\trch-\frac 2r)\nabb_A (\bk_{AN})\\ &=-\nabb_A \bigg(\bk_{AN} (\trch -\frac 2r)\bigg) +
\Theta \cdot D_*\pr H .
\endaligned
$$
Thus
\beaa
\divv \big(\ddd_4\muv &+& \half \trch\muv_A -\chih_{AB}\muv_B \big ) =   
\pr H\cdot \ddd_{4}\muv + \nabb_A\bigg(2\rr_{44}
-4\chih_{AB}\eta_B + \frac{2}r\pi - 2\bk_{AN}  (\trch-\frac 2r)\bigg) \\&-& 
\frac 2r( 3\rr_{34}  +2\rr)  
+\Theta \ric +\Theta \rr_* + \Theta\cdot D_*\pr H \nn \\
& +& \Theta\cdot\Theta\cdot \Theta 
 + \frac 1r \Theta\cdot\Theta + \frac 1{r^2}\pr H ,\nn
\eeaa
Since $\nabb r=0$ and $\curll\muv=0$, the corresponding curl equation 
takes the following final form:
\beaa
\curll \big(\ddd_4\muv + \half \trch\muv_A +\chih_{AB}\muv_B \big ) &=&
\pr H\cdot \ddd_{4}\muv - 2 \curll (\chih \cdot \muv)+\in^{AB} \nabb_A \big ((\trch-\frac 2r) 
\muv_B\big) \\&+&
\frac 1r \Theta\cdot \Theta + \Theta\cdot\Theta\cdot \Theta 
\eeaa

\subsection{ Estimate for $\nabb\trch$}\,\,\,\,
\medn

To estimate $\nabb \trch$ we commute( taking advantage of the
lemma \ref{2Comm}) the equation for $\trch $
with angular derivatives $\nabb$. Therefore, 

\begin{eqnarray}
\ddd_4\nabb\trch +\frac{3}{2}\trch \nabb\trch &=&-\nabb R_{44}
-\trch\nabb
\bk_{NN}  
- \bk_{NN}\nabb \trch  - 2\nabb\chih\cdot \chih   \nn \\ 
&-&\f12 n^{-1}\nabb n  (\half \trch^2 + \bk_{NN}\trch + \rr_{44})\nn
\end{eqnarray}

Using the transport lemma \ref{Tran}
we deduce 
\beaa
\|\nabb \trch\|_{L^2(\stu)}&\les &\frac{1}{r^2(t)}\int_u^t
 r(t')^2\bigg ( \|\nabb\rr_{44}\|_{L^2(S_{t', u})}+r(t')^{-1} \|\nabb\pr
H\|_{L^2( S_{t', u})}
+r(t')^{-2} \|\pr H\|_{L^2(S_{t', u})}\\ &+&
\|\chih\cdot\nabb\chih\|_{L^2(S_{t', u})} + r(t')^{-1} \|(\pr H)^2\|_{L^2(S_{t', u})}
+
\|\pr H  \ric(H)\|_{L^2(\sttu)}\bigg)
dt'
\eeaa
Consider the most dangerous term $\frac{1}{r^2(t)}\int_u^t
 r(t')^2  \|\nabb\rr_{44}\|_{L^2(S_{t', u})} dt'$. We estimate it with the help 
of the estimate \eqref{aso6} and find,

$$\frac{1}{r^2(t)}\int_u^t
 r(t')^2  \|\nabb\rr_{44}\|_{L^2(S_{t', u})} dt'\les\int_u^t
  \|\nabb\rr_{44}\|_{L^2(S_{t', u})} dt'\les \la^{-1-2\ep_0}
$$
Also, with the help of \eqref{aso5},
$$\frac{1}{r^2(t)}\int_u^t
 r(t')  \|\nabb\pr H\|_{L^2(S_{t', u})} dt'\les r^{-\f12}\|\nabb \pr H\|_{L^2(C_u)}
  \les r^{-\f12} \la^{-\f12 }
$$
All other terms are easier to treat.
Therefore,
\be{nabbtrchi}
r^{\half} \|\nabb\trch\|_{L^2(\stu)}\les r^{\half} \la^{-1- 2 \ep_0} +
\la^{-\f12} +  \int_u^t r(t')^{\half} \|\chih\cdot \nabb  \chih\|_{L^2(\sttu)}
dt'.
\end{equation}

It remains to  estimate $\nabb\chi$. We do this   with the help of proposition
\ref{2Hodge} applied to the Codazzi equation
\eqref{Codaz} written in the form \eqref{CZZ}. Thus
$$
\il_{S_{t,u}} |\nabb \chih |^2 + \frac{1}{r^2}|\chih|^2 \le 
 \il_{S_{t,u}} |\nabb \trch|^2+|\nabb \pr H|^2+\frac{1}{r^2}|\pr H|^2 + 
|\Theta|^4
$$
Therefore,
\begin{equation}
\begin{split}
\|\nabb \chih \|_{L^2(\stu)}&\le \|\nabb \trch \|_{L^2(\stu)}
 +\|\nabb\pr H\|_{L^2(\stu)}\\ & +\|\pr H(t)\|_{L^\infty(\stu)} +
\|\Theta\|^{2-\frac q2}_{L^\infty_x} \|\Theta\|_{L^q(\stu)}^{\frac q2}
\end{split}
\label{hodge100}
\end{equation}
for some $q>2$. Observe that we can take $q$ sufficiently close to 2 and 
use the already proved estimates \eqref{trih2} to obtain
$\|\Theta\|_{L^q(\stu)}^{\frac q2}\les \la^{-3\eps_0}$. In addition, observe 
that by H\"older inequality and \eqref{trih1}
\begin{equation}
\il_0^s \|\Theta\|^{4-q}_{L^\infty_x}\les 
s^{\frac {q-2}2}\|\Theta\|^{4-q}_{L^2_t L^\infty_x}\les \la^{-1-6\eps_0}
\label{ntr1}
\end{equation}
for all values of $q$ sufficiently close to 2.
 
Using  \eqref{hodge100}
we estimate,
\beaa
\int_u^tr(t')^{\half}
\|\chih\cdot \nabb  \chih\|_{L^2(\sttu)}dt'&\les&\int_u^t
r(t')^{\f12}\|\chih\|_{L^\infty(\sttu)}
\|\nabb \chih\|_{L^2(\sttu)}dt'\\
&\les& \int_u^t r(t')^{\half} \|\chih\|_{L^\infty(\sttu)}
\|\nabb \trch\|_{L^2(\sttu)}dt'\\
 &+& r^{\half}(t) \|\chih\|_{L^2_t L^\infty_x }
\bigg(\|\nabb \pr H\|_{L^2(C_u)}+\|\pr H\|_{L^2_t L^\infty_x} + 
\la^{-\f12 -4\eps_0}\bigg)\\
&\les& \int_u^t r(t')^{\half} \|\chih\|_{L^\infty(\sttu)}
\|\nabb \trch\|_{L^2(\sttu)}\,dt' + r^{\half}(t) \la^{-1-4\eps_0}
\eeaa
Here we have used \eqref{aso1}, \eqref{aso5}, \eqref{trih1}, \eqref{ntr1},
and the fact that $r(t')\le c r(t)$ for all $t'\le t$, which follows 
from the comparison $r(t')\approx t'-u$ and the monotonicity of
$t'-u$ along the cone $C_u$.
Therefore, returning to \eqref{nabbtrchi}, we obtain,
$$
r^{\half} \|\nabb\trch\|_{L^2(\stu)}\les r^{\half} \la^{-1-2\ep_0} +
\la^{-\frac 12} + \int_u^t r(t')^{\half} \|\chih\|_{L^\infty(\sttu)}
\|\nabb \trch\|_{L^2(\sttu)}\,dt'. 
$$
Thus, by Gronwall inequality, and the fact
that $\int_u^t\|\chih\|_{L^\infty(\sttu)}\,dt'\les 
\|\chih\|_{L_t^1L_x^\infty}\les \la^{-3\ep_0}$, we infer that
\be{nabbtrchi2}
r^{\half} \|\nabb\trch\|_{L^2(\stu)}\les r^{\half} \la^{-1- 2\ep_0} +
\la^{-\frac 12}. 
\end{equation}
Consequently, since the time interval $[0,t_*]$ obeys
$t_*\le \la^{1-8\eps_0}$, we have 
\begin{equation}
\|\sup_{r(t)\ge \frac t2}\|\nabb \trch\|_{L^2(\stu)}\|_{L^1_t}\le \la^{-3\eps_0 }
\label{ipb1}
\end{equation}
This  establishes the first part of the estimate 
\eqref{2trih8}.

\subsection{ Estimates for  $\Lb(\trch)$}
\medn

Recall the relation between $\Lb(\trch)$ and $\mu$:
$$
\mu=\Lb(\trch)-\f12 (\trch)^2 - 
\big (k_{NN}+n^{-1}N(n)\big )\trch
$$
Observe also that $\Lb (r) = \frac 1{8\pi r} \il_{\stu} \tr\chib$. Thus
$$
\Lb\big (\frac 2r\big) = -\frac 1{4\pi r^3}   \il_{\stu} \tr\chib = \frac 2{r^2} +
\frac 1{4\pi r^3} \il_{\stu} (\trch -\frac 2r +2\tr k) =  \frac 2{r^2} + 
\frac 1r\Theta
$$
In addition, $|\frac 12 (\trch)^2 -  \frac 2{r^2}|\les \frac 1r \Theta$. 
Therefore,
\bea
\|\Lb(\trch -\frac 2r)\|_{L^2(\stu)} &\les& \|\mu\|_{L^2(\stu)} + 
\|\pr H\cdot \Theta \|_{L^2(\stu)} + \frac 1r \|\Theta\|_{L^2(\stu)}\nn\\
&\les & \|\Theta\|_{L^\infty(\stu)} +  \|\mu\|_{L^2(\stu)} 
\label{100et6}
\end{eqnarray}
Here we have used  the H\"older inequality combined with the
estimate \eqref{aso2}:
$$
\|\pr H\|_{L^2(\stu)}\les 
\la^{-4\eps_0}.
$$

It remains to estimate $\|\mu\|_{L^2(\stu)}$.
We obtain this estimate from 
the transport equation \eqref{D4tmu} for $\mu $ which combined with 
Corollary \ref{corA4B3} 
can be written in the form:
\bea
L(\mu+\rr_{44}) + \trch (\mu +\rr_{44})&=& \Theta \nabb (\trch)  + \Theta \nabb\eta
+2 N(\rr_{44}) + \Theta^3 + \frac 1r\Theta^2 + \frac 1{r^2} \pr H\nn\\
&+&\frac 1r\ric (H) + \Theta \ric (H) + \Theta \rr_* + \frac 1r \rr_*\label{8.120}
\eea
\begin{remark} 
In the derivation of \eqref{8.120} we have expressed $\Lb(\rr_{44})$
in the form $L(\rr_{44})-2N(\rr_{44})$.
\end{remark}
Using the transport lemma and the estimate \eqref{aso6},
$\int_u^t\|\nab\rr_{44}\|_{L^2(\sttu)}\,dt'\les \la^{-1 - 2\eps_0}$, 
we infer that,
\begin{equation}
\begin{split}
\|\mu\|_{L^2(\stu)} &\les  \|\rr_{44}\|_{L^2(\stu)} +
\frac 1{r(t)} \il_u^t r(t') \|\Theta\|_{L^\infty(\sttu)} 
\|\nabb \eta \|_{L^2(\sttu)}\,dt' \\ 
& + \frac 1{r(t)^a}
\|\Theta\|_{L^2_t L^\infty_x} \|r(t')^a \nabb (\trch)\|_{L^2(C_u)} + 
\la^{-1-2 \eps_0}\\ &+   
\|\Theta\|^2_{L^2_t L^\infty_x}\|\Theta\|_{L^2(\stu)}  +  
\|\Theta\|^2_{L^2_t L^\infty_x} + r(t)^{-\half}  \|\pr H\|_{L^2_t L^\infty_x}\\
& + \|\ric (H)\|_{L^1_tL^\infty_x} + 
\|\Theta\|_{L^2(\stu)} \|\ric (H) \|_{L^1_t L^\infty_x} \\ &+ 
\|\Theta\|_{L^2_t L^\infty_x} \|\rr_* \|_{L^2(C_u)} + r(t)^{-\f12}\|\rr_*\|_{L^2(C_u)} 
\\ &\les  
\frac 1{r(t)} \il_u^t r(t') \|\Theta\|_{L^\infty(\sttu)} \|\nabb \eta \|_{L^2(\sttu)}\,dt'
+ \la^{-1} + r(t)^{-\half} \la^{-\f12}
\end{split}
\label{10et5}
\end{equation}
Here we have repeatedly used the H\"older inequality, 
the assumptions on the metric \eqref{aso1}-\eqref{aso7},
the already proved estimates \eqref{trih1}-\eqref{trih2} for 
$\Theta$, and the estimate\footnote{Constant $a$ can be chosen arbitrarily
from the interval $(0,2)$. Its only purpose is to remove the logarithmic
divergence at $\rho=0$.}
$$
\aligned 
\|r(t')^a \nabb (\trch)\|_{L^2(C_u)}& = 
\bigg(\il_u^t \|r(t')^a \nabb\trch\|^2_{L^2(\sttu)}\,dt'\bigg)^{\half}\\ &\les
\bigg(\il_u^t r(t')^{2a} (\la^{-1-2\eps_0} + r(t')^{-\half} \la^{-\half})^2\,
dt'\bigg)^{\half}
\\ &\les r(t)^{\half + a} \la^{-1-2\eps_0} +  r(t)^{a} \la^{-\half}\les r(t)^{a} \la^{-\half}
\endaligned
$$
following from the estimate for $\|\nabb\trch\|_{L^2(\stu)}$ proved in 
\eqref{nabbtrchi2}.

On the other hand, $\eta$ is the solution of
the Hodge system \eqref{diveta}--\eqref{curleta}:
\beaa
\divv\,\eta &=& \half\bigg(\mu +2\bk_{NN}\trch   -2|\eta|^2 -|\chih|^2
-2k_{AB}\chi_{AB}\bigg)
 - \half \de^{AB}\rr_{A{43}B},\\
\curll\,\eta &=& \half\in^{AB} k_{AC} \chih_{CB} -
 \half \in^{AB}\rr_{A{43}B}.
\eeaa
The elliptic estimate of proposition \ref{2Hodge} applied to this div-curl system
gives us the bound
\begin{equation}
\begin{split}
\|\nabb\eta\|_{L^2(\stu)} + \frac 1r \|\eta\|_{L^2(\stu)} &\les
\|\mu\|_{L^2(\stu)} + \|\Theta\|_{L^\infty(\stu)} \|\Theta\|_{L^2(\stu)}\\
& + \|\Theta\|_{L^\infty(\stu)} + \|\rr_{A{43}B}\|_{L^2(\stu)}
\end{split}
\label{10et4}
\end{equation}
Recall that according to \eqref{aso7}, 
$\|\rr_{A{43}B}\|_{L^2(C_u)}\les \la^{-\f12}$.
Thus substituting estimate \eqref{10et4} into \eqref{10et5} we obtain
$$
\|\mu\|_{L^2(\stu)}\les  \frac 1{r(t)} \il_u^t r(t') \|\Theta\|_{L^\infty(\sttu)} 
\|\mu \|_{L^2(\sttu)}\,dt' + 
\la^{-1} + r(t)^{-\half} \la^{-\half}
$$
We rewrite the above inequality in a more convenient form:
$$
r(t)^{\half} \|\mu\|_{L^2(\stu)}\les  \il_u^t  \|\Theta\|_{L^\infty(\sttu)} 
r(t')^{\half} \,\|\mu \|_{L^2(\sttu)}\,dt'+ 
r(t)^{\half} \la^{-1} + \la^{-\half}
$$
Since 
$$
\il_u^t  \|\Theta\|_{L^\infty(\sttu)}\,dt'\le 
\il_0^t \|\Theta\|_{L^\infty_x}\,dt \les t^{\half}
\|\Theta\|_{L^2_t L^\infty_x}\les \la^{-4\eps_0},
$$ 
application of 
Gronwall's inequality yields the estimate
$$
r(t)^{\half} \|\mu\|_{L^2(\stu)}\les r(t)^{\half} \la^{-1} + \la^{-\half}
$$
Returning to \eqref{100et6} we obtain
$$
\|\Lb(\trch -\frac 2r)\|_{L^2(\stu)} \les \la^{-1} + \|\Theta\|_{L^\infty_x} +
r^{-\half} \la^{-\half}.
$$
Similarly to \eqref{ipb1} we then derive the following
estimates in the exterior region:
\be{ipb4}
 \|\sup_{r\ge \frac t2}
\|\Lb(\trch - \frac 2r)\|_{L^2(\stu)}\|_{L^1_t}\le \la^{-3\eps_0 }
\end{equation}
This proves the first part of the estimate \eqref{2trih6}.

To finish the proof of \eqref{2trih6}-\eqref{2trih8}.
we first recall that 
$\Lb\big(\frac 2r\big) = \frac 2{r^2} + \frac 1r \Theta$.
Observe also that 
$$
\Lb \big(n(t-u)\big)=  n^{-1} \Lb (n) n(t-u) + n (n^{-1} - 2b^{-1})
= -1 + 2n(b^{-1} - n^{-1}) +  n^{-1} \Lb (n) n(t-u)
$$
According to Corollary  \ref{CoSn} 
$|b-n|\les s \MH$. Since by lemmas \ref{Sntu}, \ref{Rs}
the quantities $r, s$, and $n(t-u)$ are comparable, we infer that
$$
\aligned
\Lb \bigg (\frac 2r\bigg) - \Lb \bigg(\frac 2{n(t-u)}\bigg)&=
\frac 2{r^2} - \frac 2{n^2(t-u)^2} + \frac 1r \MH +\frac 1r\Theta\\ &=
2\big (\frac 1r + \frac 1{n(t-u)}\big) \big(\frac 1r -\frac 1{n(t-u)}\big) + 
 \frac 1r (\MH +\Theta)
\endaligned
$$  
Thus using Corollary \ref{Comp}, \eqref{aso1}, and \eqref{trih1} together
with the estimate for the maximal function we obtain
$$
\aligned
\|\sup_{r\ge \frac t2} \|\Lb\big (\frac 2r\big) -
\Lb \big (\frac 2{n(t-u)}\big)\|_{L^2(\stu)}\|_{L^2_t}& \les
\|\frac 1r -\frac 1{n(t-u)}\|_{L^2_t L^\infty_x}\\ &+
\|\MH\|_{L^2_t} + \|\Theta\|_{L^2_t L^\infty_x} \les 
\la^{-\frac 12 -4\eps_0}
\endaligned
$$
The above inequality followed by H\"older and \eqref{ipb4} allow us to
conclude that
\begin{equation}
\label{ipb5}
\|\sup_{r\ge \frac t2}
\|\Lb(\trch -\frac 2{n(t-u)})\|_{L^2(\stu)}\|_{L^1_t}\le \la^{-3\eps_0 },
\end{equation}
Similarly,
$$
\nabb \big(n(t-u)\big)=  n^{-1} \nabb (n) n(t-u)
$$
and consequently,
$$
\|\sup_{r\ge \frac t2}\|\nabb \big (\frac 2{n(t-u)}\big)\|_{L^2(\stu)}\|_{L^2_t}
\les \|\sup_{r\ge \frac t2}\frac 1r\|\pr H\|_{L^2(\stu)}\|_{L^2_t}
\les \|\pr H\|_{L^2_t L^\infty_x}\les \la^{-\f12 -4\eps_0}
$$
Thus we can complement \eqref{ipb1} with the estimate
\begin{equation}
\|\sup_{r\ge \frac t2}\|\nabb \big (\trch-\frac 2{n(t-u)}\big)
\|_{L^2(\stu)}\|_{L^1_t}\le \la^{-3\eps_0 },
\label{ipb7}
\end{equation}

\bigskip

It only remains to discuss  the weak estimates \eqref{lastz}. These
are a lot easier to prove and can be derived directly from
the transport equations for $\trch$ and $\chih$( see proposition \ref{proptransp}),
 in the case of the tangential derivatives $\nabb\trch$,  and 
from the transport equation for $\eta$( see proposition \ref{proptransp}), in the case of
$\Lb$  derivative\footnote{We can express $\Lb \trch$ in terms of
  $\nabb \eta$, see definition of $\mu$, and estimate the latter with
the help of the transport equation for $\eta$. }.

\end{proof}

\end{document}